\documentclass[10pt]{article}
\usepackage{amsmath,amssymb}

\newtheorem{theorem}{Theorem}[section]
\newtheorem{lemma}[theorem]{Lemma}
\newtheorem{proposition}[theorem]{Proposition}

\newtheorem{example}[theorem]{Example}

\newtheorem{corollary}[theorem]{Corollary}

\newcommand{\qed}{\enspace\vrule height6pt width4pt depth2pt}
\newenvironment{proof}{\par\noindent{\bf Proof.}}{$\qed$\par\bigskip}

\newcommand{\Z}{{\rm Z}}
\newcommand{\rk}{{\rm rk}}
\newcommand{\gr}{{\rm gr}}
\newcommand{\GK}{{\rm GK}}
\newcommand{\clKdim}{{\rm clKdim}}
\newcommand{\Spec}{{\rm Spec}}
\newcommand{\h}{{\rm h}}
\newcommand{\U}{{\rm U}}
\newcommand{\pl}{{\rm pl}}
\newcommand{\supp}{{\rm supp}}
\newcommand{\hht}{{\rm ht}}

\date{}
\begin{document}
\title{Primes of height one  and a class of Noetherian finitely presented
algebras\thanks{Research partially supported by the Onderzoeksraad
of Vrije Universiteit Brussel, Fonds voor Wetenschappelijk
Onderzoek (Flanders), Flemish-Polish bilateral agreement
BIL2005/VUB/06 and a MNiSW research grant N201 004 32/0088
(Poland).}}
\author{ Isabel Goffa\footnote{Research funded by a Ph.D grant of
the Institute for the Promotion of Innovation
through Science and Technology in Flanders
(IWT-Vlaanderen). \newline 2000 Mathematics
Subject Classification. Primary 16P40, 16H05,
16S36; Secondary 20M25, 16S15. } \and Eric
Jespers \and Jan Okni\'{n}ski} \maketitle
\date{}

\begin{abstract}
Constructions are given of Noetherian maximal orders that are
finitely presented algebras over a field $K$, defined by monomial
relations. In order to do this, it is shown that the underlying
homogeneous information determines the algebraic structure of the
algebra. So, it is natural to consider such algebras as semigroup
algebras $K[S]$ and to investigate the structure of the monoid
$S$. The relationship between the prime ideals of the algebra and
those of the monoid $S$ is one of the main tools. Results
analogous to fundamental facts known for the prime spectrum of
algebras graded by a finite group are obtained. This is then
applied to  characterize a large class of prime Noetherian maximal
orders that satisfy a polynomial identity, based on a special
class of submonoids of polycyclic-by-finite groups. The main
results are illustrated with new constructions of concrete classes
of finitely presented algebras of this type.
\end{abstract}

Because of the role of Noetherian algebras in the
algebraic approach in noncommutative geometry, new
concrete classes of finitely presented algebras
recently gained a lot of interest. Via the
applications to the solutions of the Yang-Baxter
equation, one such class also widens the interest
into other fields, such as mathematical physics.
These algebras are finitely generated, say by
elements $x_{1},\ldots ,x_{n}$, and they have a
presentation defined by ${n\choose 2}$  monomial
relations of the form $x_{i}x_{j}=x_{k}x_{l}$ so
that every word $x_{i}x_{j}$ appears at most once in
a relation. Such algebras have been extensively
studied, for example in
\cite{eting-gur,eting,gat-van,jes-okn-bin,rump}. It
was shown that they  have a rich algebraic structure
that resembles that of a polynomial algebra in
finitely many commuting generators; in particular,
they are prime Noetherian maximal orders. Clearly,
these algebras can be considered as semigroup
algebras $K[S]$, where $S$ is a monoid defined via a
presentation as that of the algebra. It turns out
that these algebras are closely related to group
algebras since $S$ is a submonoid of a torsion-free
abelian-by-finite group.

In this paper we further explore new constructions
of Noetherian maximal orders that are finitely
presented algebras which are defined via monomial
relations. In other words, we look for new
constructions of  semigroup algebras $K[S]$ of this
type. In general, it remain unsolved problems to
characterize when an arbitrary semigroup algebra is
Noetherian and when it is a prime Noetherian maximal
order. Recall that the former question even for
group algebras has been unresolved. The only class
of groups for which a positive answer has been given
is that of the polycyclic-by-finite groups. Hence it
is natural to consider the problems first for
semigroup algebras $K[S]$ of submonoids $S$ of
polycyclic-by-finite groups.  A structural
characterization of such algebras $K[S]$ that are
right Noetherian was obtained by the authors in
\cite{jes-okn-poly,jes-okn-acc}. Namely, $S$ has a
group of quotients $G=SS^{-1}$ in which there is a
normal subgroup $H$ of finite index so that
$[H,H]\subseteq S$ and $S\cap H$ is finitely
generated, or equivalently, $S$ has a group of
quotients $G$ that contains normal subgroups $F$ and
$N$ so that $F\subseteq N$, $F$ is a subgroup of
finite index in the unit group $\U(S)$ of $S$, the
group $N/F$ is abelian, $G/N$ is finite and $S\cap
N$ is finitely generated. It follows that $K[S]$ is
right Noetherian if and only if $K[S]$ is left
Noetherian. We simply say that $K[S]$ is Noetherian.
Also, if $S$ is a submonoid of a finitely generated
group $G$ that has an abelian subgroup $A$ of finite
index, then $K[S]$ is Noetherian if and only if
$S\cap A$ is finitely generated. Recall that if
$K[S]$ is Noetherian then $S$ is finitely generated
and $K[S\cap G_{1}]$ is Noetherian  for every
subgroup $G_{1}$ of $G$. If, furthermore, $G_{1}$ is
of finite index in $G$ then $K[S]$ is a finitely
generated right (and left) $K[S\cap G_{1}]$-module.
Due to a result of Chin and Quinn \cite{chin-quinn}
on rings graded by polycyclic-by-finite groups,
$K[S]$ is Noetherian if and only if $S$ satisfies
the ascending chain condition on right
(equivalently, left) ideals.

Concerning the second problem. If $S$ is  an
abelian monoid, Anderson \cite{and1,and} proved
that $K[S]$ is a prime Noetherian maximal order
if and only if $S$ is a submonoid of a finitely
generated torsion free abelian group, $S$
satisfies the ascending chain condition on ideals
and $S$ is a maximal order in its group of
quotients. More generally, Chouinard \cite{cho}
proved that a commutative monoid algebra $K[S]$
is a Krull domain if and only if $S$ is a
submonoid of a torsion free abelian group which
satisfies the ascending chain condition on cyclic
subgroups and $S$ is a Krull order in its group
of quotients. In particular, it turns out that
the height one prime ideals of $K[S]$ determined
by the minimal primes of $S$ are crucial. Brown
characterized group algebras $K[G]$ of a
polycyclic-by-finite group $G$ that are
Noetherian prime maximal orders
\cite{brown1,brown2}. This turns out to be always
the case if $G$ is a finitely generated torsion
free abelian-by-finite group (equivalently,
$K[G]$ is a Noetherian domain that satisfies a
polynomial identity, PI for short).
 In this situation all height one primes are
principally generated by a normal element. So, in the terminology
of Chatters and Jordan \cite{cha-jor}, $K[G]$ is a unique
factorization ring. The authors described in \cite{jes-okn-max}
when a semigroup algebra of a submonoid $S$ of a finitely
generated abelian-by-finite group is a Noetherian maximal order
that is a domain. The description is fully in terms of the monoid
$S$. Examples of such monoids are the binomial monoids and more
general monoids of $I$-type. As mentioned earlier in the
introduction, the latter were introduced by Gateva-Ivanova and Van
den Bergh in \cite{gat-van} and studied in several papers (see for
example \cite{eting-gur,eting,gat-sol,jes-okn-bin,rump}) and
generalized to monoids of $IG$-type in \cite{gof-jes}.

In this paper we continue the investigations of
prime Noetherian maximal orders $K[S]$ for
submonoids $S$  of a polycyclic-by-finite group
$G=SS^{-1}$. Of course, knowledge of prime ideals
is fundamental. The primes not intersecting $S$
come from primes in the group ring $K[SS^{-1}]$.
These are rather well understood through the work
of Roseblade (see for example \cite{passman}),
and in particular, the height one primes can be
handled via Brown's result. The crucial point in
our investigations are thus the primes of $K[S]$
intersecting $S$ non-trivially. In case $K[S]$ is
Noetherian and $G=SS^{-1}$ is torsion free, the
following information was proved in
\cite{jes-okn-poly}.
\begin{enumerate}
\item If $Q$ is a prime ideal of $S$ then $K[Q]$ is a prime ideal of $K[S]$.
\item The height one prime ideals of $K[S]$ intersecting $S$
non-trivially are precisely the ideals $K[Q]$
with $Q$ a minimal prime ideal of $S$.
\end{enumerate}
Note that $G$ being torsion free is equivalent
with $K[G]$ being a domain (see
\cite[Theorem~37.5]{pas-cro}).

In the first part of this paper we continue the
study of prime ideals in case $K[S]$ is prime (and
thus  $G$ is not necessarily torsion free). We show
that the first property does not remain valid,
however the second part on primes of height one
still remains true. As a consequence, we establish
going up and going down properties between prime
ideals of $S$ and prime ideals of $S\cap H$, where
$H$ is a subgroup of finite index in $G$. These are
the analogs of the important results (see
\cite[Theorem~17.9]{pas-cro})  known on the prime
ideal behaviour between a ring graded by a finite
group and its homogeneous component of degree $e$
(the identity of the grading group). As an
application it is shown that the classical Krull
dimension $\clKdim (K[S])$ of $K[S]$ is the sum of
the prime dimension of $S$ and the plinth length of
the unit group $\U(S)$. Also, a result of Schelter
is extended to the monoid $S$: the prime dimension
of $S$ is the sum of the height and depth of  any
prime ideal of $S$.

The information obtained on primes of height one
then allows us in the second section to determine
when a semigroup algebra $K[S]$ is a prime
Noetherian maximal order provided that $G=SS^{-1}$
is a finitely generated abelian-by-finite group
(that is, $K[S]$ satisfies a polynomial identity).
The result reduces the problem to the structure of
the monoid $S$ (in particular $S$ has to be a
maximal order within its group of quotients $G$) and
to that of $G$. The characterization obtained
generalizes the one given in \cite{jes-okn-max} in
case $G$ is torsion free and it shows that the
action of $G$ on minimal primes of some abelian
submonoid of $S$ is very important (as was also
discovered for the special case where $S$ is a
finitely generated monoid of $IG$-type in
\cite{gof-jes}).   Finally, in the last section, we
prove a useful criterion for verifying the maximal
order property of such $S$. We then show how this,
together with our main result, can be used to build
new examples of finitely presented algebras $R$
(defined by monomial relations) that are maximal
orders. Since the main result of the paper deals
with semigroup algebras  $K[S]$ of submonoids $S$ of
groups, we first have to check that $R$ is defined
by  such submonoids. So, in particular, in the third
section we show how to go from the language of
presentations of algebras to the language of monoids
and their semigroup algebras.

\section{Prime spectrum}

Let $K$ be a field. Recall that if $S$ is a
monoid with a group of quotients $G=SS^{-1}$ then
the semigroup algebra $K[S]$ is prime if and only
if $K[G]$ is prime
(\cite[Theorem~7.19]{okn-book1}), or equivalently
$G$ does not contain nontrivial finite normal
subgroups. The latter can also be stated as
$\Delta^{+}(G) =\{ 1 \}$, where $\Delta^{+}(G)$
is the characteristic subgroup of $G$ consisting
of the periodic elements with finitely many
conjugates (see \cite[Theorem 5.5]{pas-cro}). By
$\Delta (G)$ one denotes the subgroup of $G$
consisting of the elements with finitely many
conjugates.

Also recall that an equivalence relation $\rho$
on a semigroup $S$ is said to be a congruence if
$s\rho t$ implies $su\rho tu$ and $us\rho ut$ for
every $u\in S$. The set of $\rho$-classes,
denoted $S/\rho$, has a natural semigroup
structure inherited from $S$. In case $S$ has a
group of quotients $G$ and $H$ is a normal
subgroup of $G$ we denote by $\sim_{H}$ the
congruence on $S$ defined by: $s\sim_{H}t$ if and
only if $s=ht$ for some $h\in H$. If $G$ is a
polycyclic-by-finite group then by $\h (G)$ we
mean the Hirsch rank of $G$.

We often will make use of the following fact (\cite[Lemma~3.1]{search}). If $S$ is a submonoid of
a finitely generated abelian-by-finite group $G$ then $S$ has a group of quotients $S\Z(S)^{-1}=\{
sz^{-1} \mid s\in S,\; z\in \Z(S)\}$. It follows that if $T$ is a submonoid of a group $G$ that
contains a normal subgroup $F$ so that $F\subseteq T$ and $G/F$ is  finitely generated and
abelian-by-finite (for example $T\subseteq G$, $G$ is polycyclic-by-finite and $T$ satisfies the
ascending chain condition on right ideals, see \cite{jes-okn-acc})  then $T$ has a group of
quotients that is obtained by inverting the elements $u$ of $T$ that are central modulo $F$. In
particular, every (nonempty) ideal of $T$ contains such an element $u$.

\begin{theorem} \label{thm-primes}
Let $S$ be a submonoid of a polycyclic-by-finite
group, say with a group of quotients $G$, and let
$K$ be a field. Assume that $S$ satisfies the
ascending chain condition on right ideals and $G$
does not contain nontrivial finite normal
subgroups. Then the height one prime ideals $P$
of $K[S]$ such that $P\cap S \neq \emptyset$ are
exactly the ideals of the form $P=K[Q]$, where
$Q$ is a minimal prime ideal of $S$.
\end{theorem}
\begin{proof}
Let $P$ be a height one prime ideal of $K[S]$ such that $Q=S\cap P \neq \emptyset$. By
\cite{okn-prime} (or see \cite[Theorem~4.5.2]{jes-okn-book}, and also
\cite[Corollary~4.5.7]{jes-okn-book}), $S/Q$ embeds into ${\mathcal M}_{n}(H)$, the semigroup of
$n\times n$ monomial matrices over a group $H=TT^{-1}$ for a subsemigroup $T$ of $S$ such that
$T\cap Q =\emptyset$. Furthermore, $K[S]/K[Q]$ embeds in the matrix algebra $M_{n}(K[H])$ and the
latter is a localization of $K[S]/K[Q]$ with respect to an Ore set that does not intersect
$P/K[Q]$. Subsemigroups of $S$ not intersecting $Q$ can be identified with subsemigroups of
${\mathcal M}_{n}(H)$. It is also known that there exists an ideal $I$ of $S$ containing $Q$ such
that the nonzero elements of $I/Q$ are the matrices in $S/Q\subseteq {\mathcal M}_{n}(H)$ with
exactly one nonzero entry. Furthermore, $T$ may be chosen so that $T\subseteq I\setminus Q$ and
actually $T$ may be identified with $(e_{11}{\mathcal M}_{n}(H)e_{11}\cap (S/Q))\setminus \{0\}$,
where $e_{11}$ is a diagonal idempotent of rank one in ${\mathcal M}_{n}(H)$. Clearly $I/Q$ is an
essential ideal in $S/Q$. We know that $S$ contains a normal subgroup $F$ of $G$ such that $G/F$
is abelian-by-finite. Then $H\subseteq FH$ and clearly $FH\cap Q=\emptyset$ and $FH=FT(FT)^{-1}$.
Clearly $FTTF\subseteq FHF=FH$, whence $(FT)(TF)\cap Q=\emptyset$. Since $T$ is a diagonal
component of $I$ and $FT,TF\subseteq I$, by the matrix pattern on the matrices of rank one in
${\mathcal M}_{n}(H)$ it follows that $TF=FT=T$. Hence $H=\gr (T)=\gr (TF)=HF$.

Notice that we have a natural homomorphism $\phi :K[S]/K[Q]\longrightarrow K[S]/P$. Since $K[S]$
is Noetherian  and as $M_{n}(K[H])$ is a localization of $K[S]/K[Q]$ with respect to an Ore set of
regular elements that does not intersect $P/K[Q]$, there exists a prime ideal $R$ in $K[H]$ such
that $M_{n}(K[H]/R)$ is a localization of $K[S]/P$, see \cite[Theorem~9.22]{goo-war}. Moreover,
since $P/K[Q]$ is a prime ideal of $K[S]/K[Q]$, we also have that $M_{n}(K[H]/R)\subseteq
M_{n}(M_{t}(D))$ for a division ring $D$ and a positive integer $t$ such that $M_{t}(D)$ is the
ring of quotients of $K[H]/R$ and $\phi$ extends to a homomorphism $\phi
':M_{n}(K[H])\longrightarrow M_{n}(K[H]/R)$. Let $\rho_{P}$ be the congruence on $S$ defined by
the condition: $(x,y)\in \rho _{P}$ if $x-y\in P$. Since the image of a nonzero entry of a matrix
$s\in {\mathcal M}_{n}(H)$ under $\phi '$ is invertible in $K[H]/R$, the rank of the matrix $\phi
' (s)$ is equal to $t$ multiplied by the number of nonzero entries of $s$. Thus, the ideal $J$ of
matrices of rank at most $t$ in $S/\rho_{P}\subseteq M_{nt}(D)$ satisfies $\phi^{-1}(J)=I$.
Moreover, $S/\rho_{P}=\phi '(S/Q)\subseteq {\mathcal M}_{n}(GL_{t}(D))$ inherits the monomial
pattern of $S/Q\subseteq {\mathcal M}_{n}(H)$. From \cite[Corollary~4.4]{jes-okn-poly} we know
that there is a right Ore subsemigroup $U$ of $I$ such that $\h (UU^{-1})=\h (G)-1$ and $U\cap
P=\emptyset $, and also the image $U'$ of $U$ in $S/\rho_{P}$ is contained in a diagonal component
of $J$ (viewed as a monomial semigroup over $GL_{t}(D)$). Then $U$ can be identified with a
subsemigroup of $S/Q$ and it is contained in a diagonal component of $S/Q$. It thus follows that
$\h (H) \geq \h(\gr(U))=\h (G)-1$.

By the remark before the theorem we know that every ideal of $S$ contains an element of $S$ that
is central modulo $F$. So, choose $z\in Q$ such that $zx\in xzF$ for every $x\in S$. Suppose that
$z^{m}\in H$ for some positive integer $m$. Then $z^{m}=ab^{-1}$ for some $a,b\in T$. Hence
$z^{m}b\in T\cap P$. As $T\cap P=\emptyset$, we obtain a contradiction. It thus follows that $\h
(\gr(z,H)) > \h (H)\geq \h (G)-1$. Therefore $\h (\gr(z,H))=\h (G)$ and, because $G$ is
polycyclic-by-finite, we get $[G:\gr (z, H)]<\infty $. Since $zs\in szF$ for every $s\in S\cap H$,
it follows that $z(S\cap H)\subseteq (S\cap H)zF=(S\cap HF)z=(S\cap H)z$. As $H=(S\cap H)(S\cap
H)^{-1}$, it follows that $zH=Hz$. Therefore $\Delta^{+}(H)$ has finitely many conjugates in $G$.
Since $\Delta^{+}(H)$ is finite, we get that $\Delta^{+}(H)\subseteq \Delta^{+} (G)=\{ 1\}$ and
thus $K[H]$ is a prime algebra. Hence, the localization $M_{n}(K[H])$ of $K[S]/K[Q]$ is prime. So
$K[Q]$ is a prime ideal in $K[S]$. As $K[S]$ is prime and $P$ is of height one, we get that
$P=K[Q]$, as desired.

To prove the converse, let $Q$ be a minimal prime
ideal of $S$. Let $P$ be an ideal of $K[S]$
maximal with respect to $S\cap P = Q$. Clearly
$P$ is a prime ideal of $K[S]$. Again by the
remark before the theorem, we know that there
exists an element $z$ in $S$ that belongs to $Q$
and that is central modulo $F$. Then $zS=Sz$, so
$z$ is a normal element of $K[S]$. Since, by
assumption $K[S]$ is a prime Noetherian algebra,
the principal ideal theorem therefore yields a
prime ideal $P'$ of $K[S]$ that is of height one
so that $z\in P'$ and $P'\subseteq P$. By the
first part of the result, $P'=K[S\cap P']$. Since
$S\cap P'$ is a prime ideal of $S$ contained in
the minimal prime ideal $S\cap P = Q$, we get
that $S\cap P' =Q$. So $K[Q]=P'$ is a height one
prime ideal of $K[S]$.
\end{proof}

Recall that the rank $\rk (S)$ of a monoid $S$
(not necessarily cancellative) is the supremum of
the ranks of the free abelian subsemigroups of
$S$. The dimension $\dim (S)$ of $S$ is defined
in the following way. By definition $\dim S = 0$
if $S=\{ e\}$. If $S$ has no zero element, then
$\dim (S)$ is the maximal length $n$ of a chain
$Q_{0}\subset Q_{1}\subset Q_{2}\subset \cdots
\subset Q_{n}$, where $Q_{0}=\emptyset$ and
$Q_{i}$ are prime ideals of $S$ for $i>0$ (note
that primes  are by definition different from
$S$), or $\infty$ if such $n$ does not exist. If
$\{ e \} \neq  S$ has a zero element, then $\dim
(S)$ is the maximal length $n$ of such a chain
with all $Q_{i}$ ($i\geq 0$) prime ideals of $S$,
or $\infty$ if such $n$ does not exist. The
spectrum $\Spec (S)$ is the set of all prime
ideals of $S$.

We give an easy example that shows that
Theorem~\ref{thm-primes}  can not be extended to
semigroup algebras that are not Noetherian. Let
$S=\{ x^{i}y^{j} \mid i> 0,\; j\in {\mathbb Z}\}
\cup \{ 1 \} $, a submonoid of the free abelian
group $\gr (x,y)$ of rank two. It is easy to see
that for a given $i>0$ and $j\in {\mathbb Z}$ the
ideal $ x^{i}y^{j}S$ contains the set $\{
x^{k}y^{j} \mid k> i,\; j\in {\mathbb Z}\} $.
Therefore $P=S\setminus \{ 1 \}$ is nil modulo
the ideal $x^{i}y^{j}S$. It follows that
$S\setminus \{ 1\}$ is the only prime ideal of
$S$. Moreover $\rk (S/P)=0< \rk (S) -1$, while
$P$ is a minimal prime ideal of $S$. Clearly
$K[P]$ is a prime ideal of $K[S]$. Let $P'$ be
the $K$-linear span of the set consisting of all
elements of the form $x^{i}(y^{j}-y^{k})$ with
$i>0$ and $j,k\in {\mathbb Z}$. Then $P'$ is an
ideal of $K[S]$ and $K[S]/P'$ is isomorphic to
the polynomial algebra $K[x]$. Therefore $K[P]$
is a prime ideal of height two (note that
$\clKdim (K[S])=\rk (S)=2$ ).

The latter equalities also follow from the
following result, which will be needed later in
the paper. Let $T$  be a cancellative semigroup
and $K$ a field. If $K[T]$ satisfies a polynomial
identity then $\clKdim (K[T])=\rk (T)=\GK(
K[T])$, \cite[Theorem~23.4]{okn-book1}. (By $\GK
(R)$ we denote the Gelfand-Kirillov dimension of
a $K$-algebra $R$.)

On the other hand, for Noetherian semigroup
algebras one can also give an example (see
Example~\ref{mainexample}) showing that
Theorem~\ref{thm-primes} cannot be extended to
prime ideals of height exceeding one.

In order to give some applications to the
behaviour of prime ideals of $S$ and those of
$S\cap H$ (with $H$ a normal subgroup of finite
index in $SS^{-1}$) we prove the following
technical lemma.

\begin{lemma} \label{invariant}
Let $S$ be a submonoid of a polycyclic-by-finite
group, say with a group of quotients $G$. Assume
that $S$ satisfies the ascending chain condition
on right ideals. Then $G$ has a poly-(infinite
cyclic) normal subgroup $H$ of finite index so
that $S\cap H$ is $G$-invariant and $H/\U(S\cap
H)$ is abelian. If $G$ is abelian-by-finite then
$H$ can be chosen to be an abelian subgroup.
\end{lemma}
\begin{proof}
It is well known that $G$ contains a
characteristic subgroup $C$ that is
poly-(infinite cyclic). Since $S$ satisfies the
ascending chain condition on right ideals,   $G$
contains a normal subgroup $F$ so that $G/F$ is
abelian-by-finite and $F\subseteq S$. Hence
$C\cap F$ is a normal subgroup of $G$ that is
poly-(infinite cyclic), $C\cap F \subseteq \U(S)$
and $G/(C\cap F)$ is abelian-by-finite. It is
thus sufficient to prove the result for the
monoid $S/\sim_{C\cap F}$ and its group of
quotients $G/(C\cap F)$. In other words we may
assume that $G$ is abelian-by-finite.

So, let $A$ be a torsion free abelian and normal
subgroup of finite index in $G$. Then $S\cap A$
is a finitely generated abelian monoid and its
group of quotients is of finite index in $G$.
Since $S$ satisfies the ascending chain condition
on right ideals, we also know \cite{jes-okn-poly}
that for every $g\in S$ and $s\in S$ there exists
a positive integer $k$ so that $gs^{k}g^{-1}\in
S$. Because $G=S\Z(S)^{-1}$ (again by the remark
before Theorem~\ref{thm-primes}), the latter
property holds for all $g\in G$. As $S\cap A$ is
finitely generated abelian and $G/A$ is finite,
it follows that there exists a positive integer
$n$ so that
   $$g(S\cap A)^{(n)}g^{-1}\subseteq S\cap A,$$
for all $g\in G$,  where by definition $(S\cap
A)^{(n)}=\{ s^{n} \mid s\in S\cap A\}$. Hence
   $$ T=\bigcap_{g\in G} g(S\cap A)^{(n)}g^{-1}$$
is a $G$-invariant submonoid of $S$. Because $A$
is the group of quotients of $S\cap A$, we  get
that each $g(S\cap A)^{(n)}g^{-1}$ has a group of
quotients that is of finite index in $G$. Since
there are only finitely many such conjugates, it
is clear that $TT^{-1}$ is of finite index in
$G$. Hence, the result follows.
\end{proof}

The notion of height $\hht (P)$ of a prime ideal
of a ring $R$ is well known. We now define  the
height $\hht (Q)$ of a prime ideal $Q$ of a
monoid $S$. If $S$ does not have a zero element
then $\hht (Q)$ is the maximal length $n$ of a
chain $Q_{0}\subset Q_{1}\subset Q_{2}\subset
\cdots \subset Q_{n}=Q$, where $Q_{0}=\emptyset$
and $Q_{1},\ldots, Q_{n-1}$ are prime ideals of
$S$. On the other hand, if $S$ has a zero
element, then $\hht (Q)$ is the maximal length of
such a chain with all $Q_{i}$ prime ideals of
$S$, $i\geq 0$. If such $n$ does not exist then
we say that the height of $P$ is infinite.

\begin{corollary} \label{primes-prop}
Let $S$ be a submonoid of a polycyclic-by-finite
group, say with a group of quotients $G$, and let
$K$ be a field. Assume that $S$ satisfies the
ascending chain condition on right ideals. If $H$
is a torsion free normal subgroup of finite index
in $G$, then
\begin{enumerate}
\item If $P$ is a prime ideal of $S$ then $P\cap H =
Q_{1}\cap \cdots \cap Q_{n}$, where $Q_{1},\ldots
, Q_{n}$ are all the primes of $S\cap H$ that are
minimal over $P\cap H$. Furthermore,  $P=M\cap S$
for any prime ideal $M$ of $K[S]$ that is minimal
over $K[P]$, $\hht (M)=\hht (K[Q_{1}]) =\cdots =
\hht (K[Q_{n}])$ and $P={\mathcal B}(S(Q_{1}\cap
\cdots \cap Q_{n})S)$ (the prime radical of
$S(Q_{1}\cap \cdots \cap Q_{n})S$). If,
furthermore, $H\cap S$ is $G$-invariant then
$Q_{i}=Q_{1}^{g}$ for some $g\in G$, $1\leq i
\leq n$.
\item  If $S\cap H$ is $G$-invariant and
$Q=Q_{1}$ is a prime ideal of $S\cap H$ then
there exists a prime ideal $P$ of $S$ so that
$P\cap H=Q_{1}\cap Q_{2}\cap \cdots \cap Q_{m}$,
where $Q_{1},\ldots , Q_{m}$ are all the  prime
ideals of $S\cap H$ that are minimal over $P\cap
H$. One says that $P$ lies over $Q$. Moreover,
each $Q_{i}=Q^{g_{i}}$ for some $g_{i}\in G$.
\item {\bf Incomparability}
Suppose $S\cap H$ is $G$-invariant, $Q_{1}$ and
$Q_{2}$ are prime ideals of $S\cap H$, and
$P_{1}$ and $P_{2}$ are prime ideals of $S$. If
$P_{1}$ lies over $Q_{1}$ and $P_{2}$ lies over
$Q_{2}$ so that $Q_{1} \subseteq Q_{2}$ and
$P_{1}\subseteq P_{2}$, then $P_{1}=P_{2}$ if and
only if $Q_{1}=Q_{2}$.
\item {\bf Going up} Assume $S\cap H$ is $G$-invariant.
Suppose $Q_{2}$ is a prime ideal of $S\cap H$ and
$P_{2}$ is a prime ideal of $S$ lying over
$Q_{2}$.
\begin{enumerate}
\item
If $Q_{1}$ is a prime ideal of $S\cap H$
containing $Q_{2}$ then there exists a prime
ideal $P_{1}$ lying over $Q_{1}$ so that
$P_{2}\subseteq P_{1}$.
\item
If $P_{1}$ is a prime ideal of $S$ containing
$P_{2}$ then there exists a prime ideal $Q_{1}$
of $S\cap H$ containing $Q_{2}$ so that $P_{1}$
lies over $Q_{1}$.
\end{enumerate}
\item {\bf Going down} Assume $S\cap H$ is $G$-invariant. Suppose
$Q_{1}$ is a prime ideal of $S\cap H$ and $P_{1}$
is a prime ideal of $S$ lying over $Q_{1}$.
\begin{enumerate}
\item
If $Q_{2}$ is a prime ideal of $S\cap H$
contained in  $Q_{1}$ then there exists a prime
ideal $P_{2}$ lying over $Q_{2}$ so that
$P_{2}\subseteq P_{1}$.
\item If $P_{2}$ is a prime ideal of $S$
contained in  $P_{1}$ then there exists a prime
ideal $Q_{2}$ of $S\cap H$ contained in $Q_{1}$
so that $P_{2}$ lies over $Q_{2}$.
\end{enumerate}
\end{enumerate}
\end{corollary}
\begin{proof}
The algebra $K[S]$ has a natural gradation by the
finite group $G/H$. Its homogeneous component of
degree $e$ (the identity of the group $G$) is the
semigroup algebra $K[S\cap H]$. Let $P$ be a
prime ideal of $S$. Let $M$ be a prime ideal of
$K[S]$ minimal over $K[P]$. Note that then
$K[S]/K[P]$ inherits a natural $G/H$-gradation,
with component of degree $e$ the algebra $K[S\cap
H]/K[P\cap H]$. Because of Theorem~17.9 in
\cite{pas-cro} on going-up and down on prime
ideals of rings graded by finite groups, one gets
that $K[S\cap H] \cap M = P_{1}\cap \cdots \cap
P_{n}$, where $P_{1},\ldots , P_{n}$ are all the
prime ideals of $K[S\cap H]$ that are minimal
over $K[P]\cap K[S\cap H]=K[P\cap H]$.
Furthermore, $\hht (P_{i})=\hht (M)$ for $1\leq
i\leq n$. Since $H$ is torsion free, we know that
$K[P_{i}\cap H]$ is a prime ideal of $K[S\cap H]$
(see the introduction). Since it clearly contains
$K[P\cap H]$, it follows that $P_{i}=K[Q_{i}]$,
with $Q_{i}=P_{i}\cap H$. Furthermore, since
$K[Q]$ is a prime ideal in $K[S\cap H]$ for every
prime ideal $Q$ of $S\cap H$, it follows that
$Q_{1}, \ldots , Q_{n}$ are all the prime ideals
of $S\cap H$ minimal over $P\cap H$. So $K[S\cap
H]\cap M =K[Q_{1}\cap \cdots \cap Q_{n}]$.
Because $K[S\cap H]/K[P\cap H]$ is Noetherian, we
know that its prime radical is nilpotent. Hence
$(Q_{1}\cap \cdots \cap Q_{n})^{k} \subseteq
P\cap H$ for some positive integer $k$. Since $H$
is of finite index in $G$, it then also follows
that $M\cap S$ is an ideal of $S$ that is nil
modulo $S(P\cap H)S$. Since $K[S]$ is Noetherian,
this yields that $(M\cap S)^{l}\subseteq S(P\cap
H)S\subseteq P$, for some positive integer $l$.
As $P$ is a prime ideal, we therefore obtain that
$M\cap S=P, P\cap H=Q_{1}\cap \cdots \cap Q_{n}$
and $P={\mathcal B}(S(Q_{1}\cap \cdots \cap
Q_{n})S)$.

Assume now that, furthermore,  $S\cap H$ is
$G$-invariant. For an ideal $I$ of $S\cap H$ put
$I^{inv}=\bigcap_{g\in G} I^{g}$, the largest
invariant ideal of $S\cap H$ contained in $I$.
Clearly $SI^{inv}=I^{inv}S$ is an ideal of $S$
and $SI^{inv}\cap (S\cap H) =I^{inv}$. It follows
that
 $$SQ_{1}^{inv} \cdots S Q_{n}^{inv}\subseteq
 P.$$
Hence $Q_{i}^{inv} \subseteq P\cap (S\cap H)=Q_{1}\cap \cdots \cap Q_{n}$  for some $i$. Because
$S\cap H$ is invariant, it follows that every $Q_{i}^{g}$ is a prime ideal of $S\cap H$ (of the
same height as $Q_{i}$). Hence, for every $1\leq j \leq n$ there exists $g\in G$ with
$Q_{i}^{g}=Q_{j}$. This proves the first part of the result.

To prove the second part, let $Q$ be a prime ideal of $S\cap H$ and suppose that $S\cap H$ is
$G$-invariant. Then, each $Q^{g}$ is a prime ideal of $S\cap H$ and $\hht (Q)=\hht (Q^{g})$. Now,
if $Q'$ is a prime ideal of $S\cap H$ containing $Q^{inv}$, then $Q^{g}\subseteq Q'$ for some
$g\in G$. Hence, if $Q'$ is a prime minimal over $Q^{inv}$ then $Q^{g}= Q'$. Clearly, for every
$h\in G$, the prime $Q^{h}$ contains a prime ideal $Q''$ of $S\cap H$ minimal over $Q^{inv}$.
Hence, by the previous, $Q''=Q^{g}\subseteq Q^{h}$, for some $g$. Since $\hht(Q^{g})=\hht(Q^{h})$,
it thus follows that $Q^{h}=Q''$. So, we have shown that the ideals $Q^{g}$ are precisely  the
prime ideals of $S\cap H$ that are minimal over $Q^{inv}$. Since $H$ is torsion free we thus get
(again see the introduction) that the ideals $K[Q^{g}]$ are precisely the prime ideals of $K[S\cap
H]$ that are minimal over $K[Q^{inv}]$. Since $SQ^{inv}=Q^{inv}S$ is an ideal of $S$, the algebra
$K[S]/K[SQ^{inv}]$ has a natural $G/H$-gradation, with $K[S\cap H]/K[Q^{inv}]$ as component of
degree $e$. Hence, by \cite[Theorem~17.9]{pas-cro}, there exists a prime ideal $M$ of $K[S]$ that
is minimal over $K[SQ^{inv}]$ and $M\cap K[S\cap H]=K[Q]\cap P_{2} \cap \cdots \cap P_{m}$, where
$K[Q], P_{2},\ldots , P_{m}$ are all the prime ideals that are minimal over $M\cap K[S\cap H]$
(and these ideals are of the same height as $M$). Hence each $P_{i}$ is minimal over $K[Q^{inv}]$
and thus is of the form $K[Q^{g}]$ for some $g$. It follows that $P=M\cap S$ is a prime ideal of
$S$ so that $P\cap H$ is an intersection of $Q$ and some
 of the prime ideals $Q^{g}$, with $g\in
G$. This proves part two.

Parts (3), (4) and (5) are now immediate
consequences of parts (1) and (2) and of the
corresponding results on going up and down for
rings graded by finite groups (see for example
\cite[Theorem~17.9]{pas-cro}).
\end{proof}

A fundamental result  (see
\cite[Theorem~19.6]{pas-cro}) on  the group
algebra $K[G]$ of a polycyclic-by-finite group
says that $\clKdim (K[G]) =\pl (G)$. For the
definition on the plinth length $\pl (G)$ of $G$
we refer the reader also to \cite{pas-cro}. In
the following result we determine relations
between the considered invariants of $S$ and
$K[S]$ for Noetherian semigroup algebras $K[S]$
of submonoids $S$ of polycyclic-by-finite groups.

\begin{corollary} \label{primes-form}
Let $S$ be a submonoid of a polycyclic-by-finite
group and let $K$ be a field. Assume that $S$
satisfies the ascending chain condition on right
ideals and let $G$ be its group of quotients. Let
$F$  be a normal subgroup of $G$ so that
$F\subseteq S$, $\U(S)/F$ is finite and $G/F$ is
abelian-by-finite. The following properties hold.
\begin{enumerate}
\item $\dim (S) = \rk (G/F)$.
\item
 $\dim (S)
=\dim (S\cap W)$ for any normal subgroup $W$ of
$G$ of finite index.
\item $\Spec (S)$ is finite.
\item  $ \clKdim(K[S])=\dim (S)+\pl
(\U(S))$.
\end{enumerate}
\end{corollary}
\begin{proof}
Because of Lemma~\ref{invariant}, the group  $G$
has a poly-(infinite cyclic) normal subgroup $H$
of finite index so that $S\cap H$ is
$G$-invariant and $H/F$ is abelian for some
normal subgroup $F$ of $H$ that is contained in
$S$. Parts (1-5) of Corollary~\ref{primes-prop}
easily yield that $\dim (S)=\dim (S\cap H)$. It
also easily is verified that $\dim (S\cap H)
=\dim (T)$ where $T=(S\cap H) /\sim_{F}$. Clearly
$TT^{-1}$ is abelian and $TT^{-1}=H/F$ because
$[G:H]<\infty$. Since $\U(S)/F$ is finite, we
obtain that $\U(T)$ is finite. Hence
 $\dim (T)=\rk (T)$ (see \cite{jes-okn-poly}).  So part (1) follows.

To prove part (2) let $W$ be a normal subgroup of
finite index in $G$. Then, $S\cap W$ inherits the
assumptions on $S$ and $W$ is the group of
quotients of $S\cap W$. So, from part (1) we
obtain that $\dim (S)=\rk (G/F)=\rk (W/(F\cap
W))= \dim (S\cap W)$ and thus $\dim (S)=\dim
(S\cap W)$.

It easily is seen that there is a natural
bijective map between $\Spec (S)$ and
$\Spec(S/\sim_{F})$. Hence to prove part (3) we
may assume that $SS^{-1}$ is abelian-by-finite.
Because of Lemma~\ref{invariant} and part (1) of
Corollary~\ref{primes-prop}, we obtain that
$\Spec(S)$ is finite if the corresponding
property holds for finitely generated abelian
monoids $A=\langle a_{1},\ldots, a_{n}\rangle$.
This is obviously satisfied. Indeed, if $P\in
\Spec(A)$ then $A\setminus P$ is a submonoid of
$A$ and $a_{i}\in P$ for some $i$. It follows
that $A\setminus P$ is generated by a proper
subset of $\{ a_{1},\ldots, a_{n} \}$, whence
$|\Spec(A)|<\infty$.

Finally, we prove part (4). From
\cite[Corollary~4.4]{jes-okn-poly} we know that
$$\clKdim(K[S])= \rk(G/F)+\pl (\U(S)).$$ So the
statement follows at once from part (1).
\end{proof}

In the following proposition it is shown that the
prime spectra of $S$ and
$S/\sim_{\Delta^{+}(SS^{-1})}$ can be identified.

\begin{proposition} \label{reduction}
Let $S$ be  a submonoid of a polycyclic-by-finite
group $G=SS^{-1}$ and assume that $S$ satisfies
the ascending chain condition on right ideals.
Let $H=\Delta^{+}(G)$, the maximal finite normal
subgroup of $G$. Let $\sim$ be the congruence
relation on $S$ determined by $H$, that is,
$s\sim t$ if and only if $sH=tH$. Then the map
$\Spec (S) \longrightarrow \Spec (S/ \!\sim )$,
defined by $P\mapsto P/\!\sim $, is a bijection.
\end{proposition}
\begin{proof}
Let $\varphi : S \longrightarrow S/\! \sim \
\subseteq G/H$ be the natural epimorphism. Let
$P$ be a prime ideal of $S$.

We claim that $P=\varphi^{-1}(\varphi (P))$. One
inclusion is obvious. To prove the converse, let
$x\in \varphi^{-1}(\varphi (P))$. Then there
exists $p\in P$ such that $x\sim p$. So $x=ph$
for some $h\in H$.  Let $s\in S$. Then
 $$(hsp)^{n}= h h^{sp} h^{(sp)^{2}} \cdots h^{(sp)^{n-1}} (sp)^{n},$$
for every $n\geq 1$. Because $h\in H$, there
exists a positive integer $k$ so that
$h^{(sp)^{k}}=h$. Define $n=(k-1)|H| +1$. Then
$$(hsp)^{n} = (h h^{sp} h^{(sp)^{2}} \cdots
h^{(sp)^{k-1}})^{|H|} (sp)^{n} = (sp)^{n}\in P.$$ It follows that $xS$ is nil modulo $P$. Because
$K[S]$ is Noetherian, from  \cite[Theorem~5.18]{goo-war} we get that $x\in P$. This proves the
claim.

The claim easily implies that $\varphi (P)$ is a
prime ideal of $S$ and the statement follows.
\end{proof}

We can now prove for  the semigroups under
consideration an analogue of Schelter's theorem
on prime affine algebras that satisfy a
polynomial identity; this on its turn yields the
catenary property (see for example
\cite[Theorem~13.10.12 and
Corollary~13.10.13]{rob}).

\begin{proposition} Let $S$ be a submonoid of a
polycyclic-by-finite group, say with a group of
quotients $G$. Assume that $S$ satisfies the
ascending chain condition on right ideals. Let
$P$ be a prime ideal of $S$. Then $\dim (S/P)+
\hht (P)=\dim (S)$. Furthermore, if $\U(S)$ is
finite, then $\dim (S/P)=\rk (S/P)$.
\end{proposition}
\begin{proof}
Again let $F$ be a normal subgroup of $G$ so that
$F\subseteq \U(S)$, $[\U(S):F]<\infty$, $S\cap F$
is finitely generated and $G/F$ is
abelian-by-finite. Because of the natural
bijection between $\Spec (S)$ and $\Spec
(S/\sim_{F})$, we may replace $S$ by $S/\sim_{F}$
and thus we may assume that $G$ is
abelian-by-finite and $\U (S)$ is finite. Hence,
by Lemma~\ref{invariant}, $G$ contains a normal
torsion free abelian subgroup $A$ so that
$[G:A]<\infty$, $T=S\cap A$ is finitely generated
and $G$-invariant.

We now first prove  that $\dim (T/Q)+ \hht
(Q)=\dim (T)$ for a prime ideal $Q$ in the
abelian monoid $T$. Since $A$ is torsion free, we
know (see the introduction) that $K[Q]$ is a
prime ideal of $K[T]$ (which is of height one if
$Q$ is a minimal prime of $T$). Hence (by
Schelter's result for finitely generated
commutative algebras) $\clKdim (K[T])=\clKdim
(K[T]/K[Q]) + \hht (K[Q])$. Note that $T\setminus
Q$ is a submonoid of $A$ and $\clKdim (K[T]/K[Q])
=\clKdim (K[T\setminus Q])$ (we will several
times use this fact without specific reference).
Consequently, by Corollary~\ref{primes-form},
$\dim (T)=\dim (T\setminus Q) + \hht (K[Q])$.
Since $\dim (T\setminus Q) =\dim (T/Q)$, we need
to prove that $\hht (K[Q]) =\hht (Q)$. We prove
this by induction on $\hht (Q)$ (note that by
Corollary~\ref{primes-form} the prime spectrum of
$T$ is finite and hence every prime contains a
minimal prime ideal of $T$). If $\hht (Q)=1$ then
the statement holds. So, assume $\hht (Q)>1$. Let
$Q_{1}$ be a prime of height one contained in
$Q$. Then, by the induction hypothesis, $\hht
(Q/Q_{1}) =\hht (K[Q/Q_{1}])$. Since
$K[(T/Q_{1})/(Q/Q_{1})]\cong K[T/Q]$, we thus get
from  Schelter's result that  $\clKdim (K[T]) -1=
\clKdim (K[T/Q_{1}]) =\clKdim (K[T/Q]) + \hht
(Q/Q_{1}) =\clKdim (K[T]) -\hht (K[Q]) +\hht
(Q/Q_{1})$. Hence $\hht (K[Q]) = \hht (Q/Q_{1})
+1 \leq \hht (Q)$. Since $K[M]$ is prime in
$K[T]$ if $M$ is prime in $T$, it is clear that
$\hht (Q)\leq \hht (K[Q])$. Hence we obtain that
$\hht (K[Q]) =\hht (Q)$, as desired.

Now let $P$ be a prime ideal of $S$. Because of
Corollary~\ref{primes-prop}, $\dim (S)=\dim
(S\cap A)$, $\dim (S/P)=\dim (T/Q)$ and $\hht
(P)=\hht (Q)$, where $Q$ is a prime ideal of $T$
and $P$ lies over $Q$. From the previous it thus
follows that
 $$\dim (S/P)+\hht (P)=\dim (S).$$

So, only the last part of the statement of the
result remains to be proven. Let $Q_{1}=Q$ and
write $P\cap T=Q_{1}\cap \cdots \cap Q_{n}$,
where $Q_{1},\ldots , Q_{n}$ are all primes
minimal over $P\cap T$. We know that $\clKdim
(K[T/Q_{i}]) =\rk (T/Q_{i})$. Furthermore,
because $TT^{-1}$ is torsion free and $K[T]$ is
Noetherian, we also know that $K[Q_{i}]$ is a
prime ideal of $K[T]$ that is minimal over
$K[P\cap T]$ (see the introduction). Hence,
$K[Q_{i}/(P\cap T)]$ is a minimal prime ideal in
the finitely generated commutative algebra
$K[T/(P\cap T)]$. It follows that
 \begin{eqnarray}
 \clKdim (K[T/(P\cap T)]) &=&
         \clKdim (K[T/Q_{i}]) + \hht (K[Q_{i}/(P\cap
         T)]) \nonumber\\
  &=&\clKdim (K[T/Q_{i}]) =\rk (T/Q_{i}).
  \label{equalb}
 \end{eqnarray}
Since $K[T/(P\cap T)]$ is a finitely generated
commutative algebra, it is well known that
$\clKdim (K[ T/(P\cap T)]) =\GK (K[T/(P\cap T)])$
(see for example \cite[Theorem~4.5]{krause}).
  From \cite[Theorem~23.14]{okn-book1}) it follows
that
 $$\rk (T/(P\cap T))= \clKdim (K[T/(P\cap T)]).$$
Using (\ref{equalb}) we thus get
 \begin{eqnarray*}
 \rk (T/(P\cap T))&=&\clKdim (K[T]/(K[Q_{1}]\cap \cdots \cap K[Q_{n}]))\\
   &=& \sup \{ \clKdim (K[T/Q_{i}])\mid 1\leq i \leq
   n\}\\
   &=&\clKdim (K[T/Q]) =\rk (T/Q).
 \end{eqnarray*}
Since $TT^{-1}$ is of finite index in $SS^{-1}$,
it follows that
 $$\rk (S/P) =\rk (T/(T\cap P)) =\rk
(T/Q).$$ Because $\U(T/Q)$ is finite,
Corollary~\ref{primes-form} yields that $\rk
(T/Q)= \dim (T/Q)$ and thus we get that
 \begin{eqnarray*}
 \rk(S/P) =\rk
(T/Q) &=& \dim (T/Q)
  =\dim (T)-\hht (Q)\\
   &=&\dim (S)-\hht (P) = \dim (S/P).
 \end{eqnarray*} This finishes the proof.
\end{proof}

\section{Maximal Orders}

In this section we describe when a semigroup
algebra $K[S]$ of a cancellative submonoid $S$ of
a polycyclic-by-finite group $G$ is a prime
Noetherian maximal order that satisfies a
polynomial identity. In case $G$ is torsion free
such a result was obtained in \cite{jes-okn-max}
and in  case $S=G$ this was done by Brown in
\cite{brown1,brown2} (even without the
restriction that $K[S]$ has to be PI).

For completeness' sake  we  recall some notation
and terminology on (maximal) orders. We state
these in the semigroup context (see for example
\cite{gil} and \cite{wau-kru}) as these are
basically the same as in the more familiar ring
case. A cancellative monoid $S$ which has a left
and right group of quotients $G$ is called an
{\sl order}. Such a monoid $S$ is called a {\sl
maximal order} if there does not exist a
submonoid $S'$ of $G$ properly containing $S$ and
such that $aS'b\subseteq S$ for some $a,b\in G$.
For subsets $A,B$ of $G$ we define $(A:_{l}B)=\{
g\in G \; |\; gB\subseteq A\} $ and by
$(A:_{r}B)=\{ g\in G \; |\; Bg\subseteq A\}$.
Note that $S$ is a maximal order if and only if
$(I:_{l}I)=(I:_{r}I)=S$ for every {\sl fractional
ideal} $I$ of $S$. The latter means that
$SIS\subseteq I$ and $cI,Id\subseteq S$ for some
$c,d\in S$. If $S$ is a maximal order, then
$(S:_{l}I)=(S:_{r}I)$ for any fractional ideal
$I$; we simply denote this fractional ideal by
$(S:I)$ or by $I^{-1}$. Recall that then $I$  is
said to be {\sl divisorial} if $I=I^{*}$, where
$I^{*}=(S:(S:I))$. The divisorial product $I*J$
of two divisorial ideals $I$ and $J$ is defined
as $(IJ)^{*}$.  Also recall that a fractional
ideal is said to be invertible if $IJ=JI=S$ for
some fractional ideal $J$ of $S$. In this case
$J=I^{-1}$ and $I$ is a divisorial ideal.

Recall then that (see for example \cite{wau-kru})
a cancellative monoid $S$ is said to be a Krull
order if and only if $S$ is a maximal order
satisfying the ascending chain condition on
divisorial integral ideals (the latter are the
fractional ideals contained in $S$). In this case
the set $D(S)$ of divisorial fractional ideals is
a free abelian group for the $*$ operation. If
$SS^{-1}$ is abelian-by-finite then (as said
before, see \cite[Lemma~1.1]{jes-wang}) every
ideal of $S$ contains a central element and it
follows that the minimal primes of $S$ form a
free basis for $D(S)$, \cite{wau-kru}.

Similarly a prime Goldie ring $R$ is said to be a
Krull order if $R$ is a maximal order that
satisfies the ascending chain condition on
divisorial integral ideals. Although there are
several notions of noncommutative Krull orders,
for  rings satisfying a polynomial identity all
these notions are the same.

In the next theorem we collect some of the
essential properties of these orders. For details
we refer the reader to
\cite{chamarie-t,chamarie-a}. For a ring $R$ and
an Ore set $C$ of regular elements in $R$ we
denote by $R_{C}$ the classical localization of
$R$ with respect  to $C$. The classical ring of
quotients of a prime Goldie ring $R$ is denoted
by $Q_{cl}(R)$. The prime spectrum of $R$ is
denoted by $\Spec (R)$, the set of height one
prime ideals of $R$ by $X^{1}(R)$.

\begin{theorem} \label{cha-kru}
Let $R$ be a prime Krull order satisfying a
polynomial identity. Then the following
properties hold.
\begin{enumerate}
\item
The divisorial ideals form a free abelian group with basis
$X^{1}(R)$, the height one primes of $R$.
\item If $P\in X^{1}(R)$ then $P\cap \Z(R)\in X^{1}(\Z(R))$, and
furthermore, for any ideal $I$ of $R$, $I\subseteq P$ if and only
if $I\cap \Z(R) \subseteq  P\cap \Z(R)$.
\item $R=\bigcap R_{Z(R)\setminus P}$, where the intersection is taken
over all height one primes of $R$, and every
regular element $r\in R$ is invertible in almost
all (that is, except possibly finitely many)
localizations $R_{\Z(R)\setminus P}$.
Furthermore, each $R_{\Z(R)\setminus P}$ is a
left and right principal ideal ring with a unique
nonzero prime ideal.
\item For a multiplicatively closed set of ideals ${\cal M}$ of $R$,
the (localized) ring $R_{{\cal M}}=\{ q\in
Q_{cl}(R) \mid Iq \subseteq R,\; \mbox{for some }
I\in {\cal M}\}$ is a Krull order, and
 $$R_{{\cal M}} =\bigcap R_{\Z(R)\setminus P},$$
where the intersection is taken over those height
one primes $P$ for which $R_{{\cal M}} \subseteq
R_{\Z(R)\setminus P}$.
\end{enumerate}
\end{theorem}

Next we prove some necessary condition for $K[S]$
to be a prime Noetherian  maximal order that
satisfies a polynomial identity.

\begin{lemma}\label{lem31}
Let $S$ be a submonoid of an abelian-by-finite
group $G=SS^{-1}$ and let $K$ be  a field. Let
$A$ be an abelian subgroup that is  normal and of
finite index in $G$ and let $P$ be a minimal
prime ideal of $S$. The following properties
hold.
\begin{enumerate}
\item
If $S$ is a maximal order then $S\cap A$ is
$G$-invariant and $S\cap A$ is a maximal order in
its group of quotients.
\item If
$K[S]$ is a prime Noetherian maximal order, then
$S$ is a maximal order, $S\cap A$  and $P\cap A$
are $G$-invariant.
\end{enumerate}
\end{lemma}

\begin{proof}
By \cite[Lemma~2.1]{jes-okn-max}, if  $S$ is a maximal order then $S\cap A$ is $G$-invariant and
it is a maximal order in its group of quotients. Assume now that $K[S]$ is a prime Noetherian
maximal order. It is straightforward (as in the proof of Lemma~3.3 in \cite{jes-okn-max}) to
verify that then $S$ is a maximal order, and thus $S\cap A$ is $G$-invariant. Let $P$ be a minimal
prime ideal of $S$. Of course, $A\subseteq \Delta (G)$ and thus $P\cap A= (P\cap \Delta (G)) \cap
(A\cap S)$. Hence to prove that $P\cap A$ is $G$-invariant, we may assume in the remainder that
$A=\Delta (G)$. Indeed, since $K[S]$ is prime, by  \cite[Theorem~7.19]{okn-book1}, $K[G]$ is
prime. Then, from \cite[Theorem~4.2.10]{passman} it follows that $\Delta(G)$ is abelian. Now, by
Corollary~\ref{primes-prop}, $P\cap A=Q_{1}\cap \cdots \cap Q_{n}$, an intersection of minimal
primes of $A$ that are $G$-conjugate. We need to show that $\{ Q_{1},\ldots , Q_{n}\} =\{
Q_{1}^{g}\mid g\in G\}$.  Suppose the contrary, so assume $Q'$ is a minimal prime of $S\cap A$
that is conjugate to $Q_{1}$ but is different from all $Q_{i}$, for $1\leq i\leq n$. Then by
Corollary~\ref{primes-prop}, there exists a prime ideal $P'$ of $S$ so that $P'\cap A=Q'\cap
Q_{2}'\cap \cdots \cap Q_{m}'$ for some minimal primes $Q_{2}',\ldots , Q_{m}'$ of $S\cap A$.
Because of Theorem~\ref{thm-primes}, both $K[P]$ and $K[P']$ are distinct height one prime ideals
of $K[S]$. As, by assumption, $K[S]$ is a prime PI Noetherian maximal order, it follows from
Theorem~\ref{cha-kru} that there exists a central element of $K[S]$ that belongs to $K[P']$ but
not to $K[P]$. Since central elements of $K[S]$ are linear combinations of finite conjugacy class
sums of $K[G]$ we get that $P'$ contains a $G$-conjugacy class $C$ that does not belong to $P$.
Since $A=\Delta (G)$, we thus get that $C\subseteq A$. Hence, $C\subseteq Q'$. As $Q'$ is a
$G$-conjugate of $Q_{1}$ it follows that $C\subseteq Q_{1}\cap \cdots \cap Q_{n}=P\cap A$, a
contradiction.
\end{proof}

We need one more lemma in order to prove the main
theorem of this section. If $\alpha =\sum_{s\in
S} k_{s} s$ (with each $k_{s}\in K$) then we
denote by $\supp (\alpha ) =\{ s\in S \mid
k_{s}\neq 0\}$ the support of $\alpha$.

\begin{lemma}\label{lem32}
Let $S$ be a submonoid of an abelian-by-finite
group $G=SS^{-1}$.  Let  $K$ be  a field and
suppose $K[S]$ is Noetherian. Let $A$ be a normal
abelian subgroup of finite index in $G$. Assume
that $P$ is a prime ideal of $S$ so that $K[P]$
is a prime ideal of $K[S]$ and $P\cap A$ is
$G$-invariant. Also assume that $S\cap A$ is
$G$-invariant. If $J$ is an ideal of $S$ not
contained in $P$ then $J\cap A$ contains a
$G$-conjugacy class $D$ such that $D\not
\subseteq P$ (and thus clearly $D\subseteq (J\cap
A)\setminus P$).
\end{lemma}
\begin{proof}
We may assume that $J$ is a proper ideal of $S$. The prime algebra $R=K[S]/K[P]$  has a natural
$G/A$-gradation, with identity component $R_{e}=K[S\cap A]/K[P\cap A]$, a semiprime commutative
algebra. Let $L$ be a nonzero ideal of $R$. It then follows from \cite[Theorem~1.7]{coh-row} that
$r_{e}\in L\cap R_{e}$, for some regular element $r_{e}$ of $R_{e}$. Because of the assumptions,
$G$ acts by conjugation in $R_{e}$. Clearly, $r$ has only finitely many such conjugates, say
$r_{1},\ldots , r_{m}$ and $r_{1}\cdots r_{m}\in L \cap R_{e}$ is central and nonzero. So, $L$
contains a non-trivial element in  $R_{e}\cap \Z(R)$.

We apply the above to the ideal
$L=(K[J]+K[P])/K[P]$. So, let $\alpha \in K[J\cap
A]$ be such that the image $\overline{\alpha} \in
K[S]/K[P]$ is nonzero and lies in the center. So
$\alpha$ is regular modulo $K[P]$. We may assume
that $\supp (\alpha)\cap P=\emptyset$. Note also
that $1\not\in \supp (\alpha )$. Write $\alpha =
\alpha_{1}+\cdots + \alpha_{q}$, where
$\alpha_{i}$ have supports contained in different
$G$-conjugacy classes. Then $g\alpha_{1}+\cdots +
g\alpha_{q}-(\alpha_{1}g+\cdots + \alpha_{q}g)
\in K[P]$ for every $g\in S$. Clearly,
$g\alpha_{i}$ and $\alpha_{j}g$ have disjoint
supports if $i\neq j$. So we must have
$g\alpha_{i}-\alpha_{i}g\in K[P]$ for every $i$.
Then every $\alpha_{i}$ also lies in the center
modulo $K[P]$ and $\supp(\alpha _{i})$ is
contained in a $G$-conjugacy class $D_{i}$. By
the hypothesis, $D_{i}\subseteq (S\cap
A)\setminus P$. Hence, replacing $\alpha $ by
$\alpha_{1}$, we may assume that $\alpha =
\alpha_{1}$. Write $D=D_{1}$.

Because of the assumptions and
Corollary~\ref{primes-prop}, $P\cap A=Q_{1}\cap
\cdots \cap Q_{n}$ where $Q_{1},\ldots ,Q_{n}$ is
a full orbit of conjugate primes in $S\cap A$.
For $1\leq i \leq n$ let  $A_{i}=\left(
\bigcap_{j,\; j\neq i} Q_{j}\right) \setminus
Q_{i}$. Each $A_{i}$ is a subsemigroup of
$A\setminus P$. Each $g\in G$ permutes the sets
$A_{1},\ldots ,A_{n}$ (by conjugation). Let
$I=\bigcup_{i=1}^{n} A_{i} \cup (P\cap A)$. Then
$I$ is $G$-invariant, $SI=IS$ and $SI\cap A=I$.
Replacing $J$ by $J\cap SI$  we may assume that
$J\subseteq SI$ and thus also $\alpha \in K[J\cap
A]\subseteq K[SI\cap A]=K[I]$. Hence, we can
write $\alpha = \gamma_{1}+\cdots +\gamma_{n}$,
with $\supp (\gamma_{i})\subseteq A_{i}$ for
$i=1,\ldots, n$. Notice that
$\supp(\gamma_{i})\neq \emptyset$ for every $i$.
Indeed, suppose the contrary, that is, assume $
\gamma_{i}=0$ for some $i$. Then, $\alpha
A_{i}\subseteq P$, in contradiction with the
regularity of $\alpha$ modulo $K[P]$. So, indeed,
$\supp (\gamma_{i})\neq \emptyset$ for every $i$.
Let $E_{i}=D\cap A_{i}$. Since also $\supp
(\gamma_{i})\subseteq D$, we get that
$\emptyset\neq \supp (\gamma_{i})\subseteq
E_{i}$. Put $a_{i}=\prod_{x\in E_{i}}x$. Clearly
$a_{i}\in A_{i}$. It follows that $a_{1}+\cdots
+a_{n}\in K[J]$. Let $1\leq i\leq n$ and $g\in
G$. Then $g^{-1}A_{i}g=A_{j}$ for some $j$. Hence
$g^{-1}E_{i}g=E_{j}$ and thus
$g^{-1}a_{i}g=a_{j}$. So, conjugation by elements
of $G$ permutes $a_{1},\ldots, a_{n}$. Therefore
the result follows.
\end{proof}

Recall from \cite{brown1,pas-cro} that a group
$G$ is dihedral free if, for every subgroup $D$
of $G$ isomorphic to the infinite dihedral group,
the normalizer $N_{G}(D)$ of $D$ in $G$ has
infinite index.

\begin{theorem}\label{thm33}
Let $K$ be a field and let $S$ be a submonoid of a finitely generated  abelian-by-finite group.
Let $A$ be an abelian normal subgroup of finite index in $G=SS^{-1}$. The following conditions are
equivalent.
\begin{enumerate}
\item $K[S]$ is a prime Noetherian maximal order.
\item
$S$ is a maximal order that satisfies the
ascending chain condition on right ideals,
$\Delta^{+}(G)=\{ 1 \}$, $G$ is dihedral free and
for every minimal prime ideal $P$ of $S$ the set
$A\cap P$ is $G$-invariant.
\end{enumerate}
\end{theorem}

\begin{proof}
Assume that $K[S]$ is a prime Noetherian maximal
order. Then, by Theorem~\ref{cha-kru}, the
localization $K[G]=K[S\Z(S)^{-1}]$ of $K[S]$ also
is a maximal order. Because of Brown's result on
the description of group algebras of
polycyclic-by-finite groups that are maximal
orders, the latter holds if and only if
$\Delta^{+}(G)=\{ 1\}$ and $G$ is dihedral free.
Lemma~\ref{lem31} then yields that the other
conditions listed in (2) hold as well.

Conversely, assume that condition (2) holds. The assumption on
$\Delta^{+}(G)$ yields that $K[G]$ and thus $K[S]$ is prime.
Because $S$ is a maximal order, Lemma~\ref{lem31} gives that
$S\cap A$ is $G$-invariant for any abelian normal subgroup $A$ of
$G$ of finite index. To prove that $K[S]$ is a maximal order one
can follow the lines of the proof of Theorem~3.5 in
\cite{jes-okn-max}. We now give a simplified proof.

We begin by showing that  if $P$ is a minimal prime ideal of $S$ then the localized ring
$R=K[S]((\Z(K[S])\cap K[A])\setminus K[P])^{-1}$ is a maximal order. To do so, we show that $R$ is
a local ring with unique maximal ideal $RP$ and so that $RP$ is invertible and every proper
nonzero ideal of $R$ is of the form $(RP)^{n}$ for some positive integer $n$. First we show that
$RP$ is the only height one prime ideal of $R$.  Of course if $Q'$ is a height one prime ideal of
$R$ then $Q=K[S]\cap Q'$ is a height one prime ideal of $K[S]$ that does not intersect
$(\Z(K[S])\cap K[A])\setminus K[P]$. Because of Theorem~\ref{thm-primes}, either $Q=K[S\cap Q]$ or
$S\cap Q =\emptyset$. Because of Lemma~\ref{lem32}, the former implies that $S\cap Q \subseteq P$
and thus $Q=K[P]$, as desired. So assume $S\cap Q=\emptyset$, or equivalently, $A\cap Q
=\emptyset$. Because $S\cap A$ is $G$-invariant and $Q$ does not contain homogeneous elements, it
follows (see for example Lemma~7.1.4 in \cite{jes-okn-book}) that $\overline{Q}=\bigcap_{g\in G}
g^{-1}(Q\cap K[A])g$ is not contained in $K[P\cap A]$ . As $\overline{Q}$ is $G$-invariant, we get
that $K[S]\overline{Q}$ is an ideal of $K[S]$. So $(K[S]\overline{Q}+K[P])/K[P]$ is a nonzero
ideal of the Noetherian algebra $K[S]/K[P]$. This algebra has a natural $G/A$-gradation, with
component of degree $e$ the semiprime algebra $K[S\cap A]/K[P\cap A]$. By Theorem~1.7 in
\cite{coh-row}, the ring $K[S]/K[P]$ has a classical ring of quotients that is obtained by
inverting the regular elements of $K[S\cap A]/K[P\cap A]$. Hence the ideal
$(K[S]\overline{Q}+K[P])/K[P]$ contains a regular element $\overline{\alpha}$ that is contained in
$K[S\cap A]/K[P\cap A]$. Since $K[S\cap A]$ and $K[P\cap A]$ are $G$-invariant, conjugation
induces an action of $G$ on $K[S\cap A]/K[P\cap A]$. Hence the product of the finitely many
conjugates of $\overline{\alpha}$ also belongs to $K[S\cap A]/K[P\cap A]$. Since this element is
central, we thus may assume that $\overline{\alpha}$ also is central in $K[S]/K[P]$ and clearly
$\overline{\alpha}\in (\overline{Q}+K[P\cap A])/K[P\cap A]$. Write $\overline{\alpha} = \gamma +
K[P\cap A]$ for some $\gamma \in \overline{Q}$. Let $\beta$ be the product of the distinct
conjugates of $\gamma$. Then $\beta +K[P\cap A] =\gamma^{m}+K[P\cap A]$ for some positive integer
$m$. Since $K[P\cap A]$ is a semiprime ideal in $K[S\cap A]$, it follows that $\beta \not\in
K[P\cap A]$. Hence $\beta \in (Q\cap \Z(K[S]))\setminus K[P]$, a contradiction. This implies that
indeed $RP$ is the only height one prime ideal of $R$.

Because of Lemma~\ref{lem31},  $S\cap A$ is a finitely generated maximal order. Hence  we know
that $K[S\cap A]$ is a Noetherian maximal order (\cite{and1,and}) and thus it is well known (or
use Theorem~\ref{cha-kru}) that $K[S\cap A]((\Z(K[S]\cap K[A])\setminus P)^{-1}$ is a Noetherian
maximal order with only finitely many height one prime ideals.  Hence  it is a principal ideal
domain (see for example \cite{fossum}) and thus it has prime dimension one. As this is the
component of degree $e$ of the $G/A$-graded ring $R$, from \cite[Theorem~17.9]{pas-cro} it follows
that $R$ also has dimension one. Hence $RP$ is the only maximal ideal of $R$. As $R$ also is a PI
algebra, we then obtain that $RP$ is the Jacobson radical of $R$ and thus, by
\cite[Theorem~8.12]{goo-war}, $\bigcap_{n}(RP)^{n}=\{ 0\}$. Also note that $P(S:P)$ is an ideal of
$S$ that is not contained in $P$. Hence $RP(S:P)=R$ and thus the unique maximal ideal  $RP$ is an
invertible ideal of $R$. It then easily follows that every proper nonzero ideal of $R$ is of the
form $(RP)^{n}$ for some unique positive integer $n$. This proves the desired properties of $R$.

Let $\alpha =\sum_{i=1}^{n}k_{i}g_{i} \in K[G]$,
where  $0\neq k_{i} \in K$ and $g_{i}\in G$ for each
$1\leq i \leq n$ and $g_{i}\neq g_{j}$ for $i\neq j$
We now show that if $\alpha  \in K[G] \cap
\bigcap_{P} K[S]((\Z(K[S])\cap K[A])\setminus
K[P])^{-1}$ (where the intersection runs over all
minimal primes $P$ of $S$) then $\alpha \in K[S]$.
We prove this by induction on $n$. If $n=1$ then
$\alpha =k_{1}g_{1}$. For each minimal prime $P$ of
$S$ there then exists a central element $\delta$ of
$K[S]$ that belongs to $K[A]\setminus P$ so that
$\delta k_{1}g_{1}\in K[S]$. Hence there is a
$G$-conjugacy class $C(P)$ so that $C\subseteq
(S\cap A)\setminus P$ and $C(P)g_{1}\subseteq S$.
Let $C=\bigcup_{P} C(P)$. Then $SCg_{1}\subseteq S$
and $SC$ is an ideal of $S$ that is not contained in
any of the minimal prime ideals of $S$. Since $S$ is
a maximal order, it follows that $g_{1} \in S$, as
desired. Now assume $n>1$. Since $\bigcap_{P}
K[S]((\Z(K[S])\cap K[A])\setminus K[P])^{-1}$ is a
$G/A$-graded ring, the induction hypothesis yields
that we may assume that $\alpha$ is
$G/A$-homogeneous, that is, each $g_{i}g_{j}^{-1}\in
A$. Since this statement holds for any normal
abelian subgroup of $G$ of finite index and because
$A$ is residually finite, we get that $g_{i}=g_{j}$
for all $i=j$. Hence we may assume $n=1$ and thus by
the above $\alpha \in K[S]$.

So we have shown that $K[S]= K[G] \cap \bigcap_{P} K[S]((\Z(K[S])\cap K[A])\setminus K[P])^{-1}$.
Recall that by Theorem~F in \cite{brown1}, since $G$ is dihedral-free, $K[G]$ is a maximal order.
Since also each  $K[S]((\Z(K[S])\cap K[A])\setminus K[P])^{-1}$ is a  maximal order and a central
localization of $K[S]$, it follows that $K[S]$ is a maximal order. This finishes the proof.
\end{proof}

\section{Constructing examples}

Theorem~\ref{thm33} reduces the problem of determining when $K[S]$
is a prime Noetherian maximal order to the algebraic structure of
$S$. It hence provides a strong tool for constructing new classes
of such algebras. For some examples the required conditions on $S$
can easily be verified, but on the other hand, for some examples
this still requires substantial work. In this section this is
illustrated with some concrete constructions.

A first class of examples consists of algebras
defined by monoids of $I$-type (see
\cite{gat-van,jes-okn-book}). Recall that, in
particular, these are quadratic algebras $R$ with
a presentation defined by $n$ generators
$x_{1},\ldots , x_{n}$ and with ${n \choose 2}$
relations of the form $x_{i}x_{j}=x_{k} x_{l}$ so
that every word $x_{i}x_{j}$ appears at most once
in one of the defining relations. Clearly
$R=K[S]$, where $S$ is the monoid defined by the
same presentation. It turns out that $S$ has a
group of quotients $G$ that is torsion free and
has a free abelian subgroup $A$ of finite index.
Furthermore, for any minimal prime ideal $P$ of
$S$ one has that $P=Ss=sS$, $P\cap A=(S\cap A)a$,
for some $s\in S,\; a\in A$, and $P\cap A$ is
$G$-invariant. Using Theorem~\ref{thm33} we then
immediately recover the known result that $R$ is
a maximal order. The only noncommutative algebra
of such type which is generated by two elements
is $K\langle x,y \mid x^{2}=y^{2}\rangle$
(\cite{gat-van}).  A related example on three
generators that is not of this type is $K\langle
x,y,z\mid x^{2}=y^{2}=z^{2}, zx=yz,zy=xz\rangle$.
As an application of Theorem~\ref{thm33} one can
show by elementary calculations that this algebra
also is a prime Noetherian PI maximal order.

In the remainder of this section we discuss in full detail one
more construction that illustrates Theorem~\ref{thm33} but also
shows that certain assumptions in Theorem~\ref{thm-primes} are
essential. Before this, we establish a useful general method for
constructing nonabelian submonoids of abelian-by-finite groups
that are maximal orders, starting from abelian maximal orders.

\begin{proposition}\label{notesjan}
Let $A$ be an abelian normal subgroup of finite index in a group
$G$. Suppose that $B$ is a submonoid of $A$ so that $A=BB^{-1}$
and $B$ is a finitely generated maximal order. Let $S$ be a
submonoid of $G$ such that $G=SS^{-1}$ and $S\cap A=B$. Then $S$
is a maximal order that satisfies the ascending chain condition on
right ideals if and only if $S$ is maximal among all submonoids $T$
of $G$ with $T\cap A=B$.
\end{proposition}
\begin{proof}
First suppose $S$ is a maximal order that satisfies the ascending
chain condition on right ideals. Suppose that $S\subseteq
T\subseteq G$ for a submonoid $T$ of $G$ such that $S\cap
A=B=T\cap A$. By assumption, $B$ is finitely generated. Since also
$A$ is normal and of finite index in $G$, it thus follows (see the
introduction) that $K[T]$ is Noetherian and it is a finitely
generated right $K[B]$-module. Thus $T$ satisfies the ascending
chain condition on one-sided ideals and $T=\bigcup _{i=1}^{n}
t_{i}B$ for some $n\geq 1$ and $t_{i}\in T$. Since $G$ is finitely
generated and abelian-by-finite, we know that for every $i$ there
exists $z_{i}\in \Z(S)$ such that $z_{i}t_{i}\in S$. Let
$z=z_{1}\cdots z_{k}$. Then $zt_{i}\in S$ and therefore
$zT=\bigcup _{i} zt_{i}(S\cap A)\subseteq S$. Since $S$ is a
maximal order, this implies that $T=S$. So, we have shown that $S$
is maximal among all submonoids $T\subseteq SS^{-1}$ such that
$T\cap A=S\cap A$.

Conversely, assume that $S\subseteq G$ is maximal among all
submonoids $T$ of $G$ with $T\cap A=B$. As above,  because $S\cap
A$ is finitely generated, $S$ satisfies the ascending chain
condition on right ideals. Suppose that $T\subseteq G$ is a
submonoid such that $S\subseteq T$ and $gTh\subseteq S$ for some
$g,h\in G$. There exist $s,t\in S$ such that $sg,ht\in \Z(G)\cap
B$. So $sgTht\subseteq S$ and therefore $Tz\subseteq S$ for some
$z\in \Z(G)\cap B$. In particular $(T\cap A)z\subseteq S\cap A=B$.
Since $S\cap A\subseteq T\cap A$ and $S\cap A$ is a maximal order,
it follows that $T\cap A=S\cap A$. Because $S\subseteq T$, the
assumption on $S$ then implies that $T=S$. Therefore $S$ is a
maximal order.
\end{proof}

In order to illustrate the above proposition with a
concrete example, we start with the following
construction of a monoid that contains the abelian
monoid generated by $a_{1},a_{2},a_{3},a_{4}$ and
defined by the extra relation
$a_{1}a_{2}=a_{3}a_{4}$. It was shown in \cite{and}
that  the latter is a cancellative monoid that is a
maximal order.

\begin{example}   \label{abel-max}
The abelian monoid  $B=\langle
a_{1},a_{2},a_{3},a_{4},a_{5},a_{6}\rangle$
defined by the relations
$a_{1}a_{2}=a_{3}a_{4}=a_{5}a_{6}$ is a
cancellative monoid that is a maximal order (in
its torsion free group of quotients).
\end{example}
\begin{proof}
Let $F= \langle x_{1},\ldots ,x_{8} \rangle $ be
a free abelian monoid of rank $8$. Define
$$b_{1}=x_{1}x_{2}x_{3}x_{4}, \,
b_{2}=x_{5}x_{6}x_{7}x_{8}, \,
b_{3}=x_{1}x_{2}x_{5}x_{6}, $$
$$b_{4}=x_{3}x_{4}x_{7}x_{8},\,
b_{5}=x_{1}x_{3}x_{5}x_{7}, \,
b_{6}=x_{2}x_{4}x_{6}x_{8}.$$ Clearly, these
$b_{i}$ satisfy the defining relations for $B$.
We claim that actually $B\cong \langle
b_{1},\ldots, b_{6}\rangle$ under the map
determined by $a_{i}\mapsto b_{i}$, $i=1,\ldots
,6$. In order to prove this, suppose that there
is a relation $w=v$, where $w,v$ are nontrivial
words in $b_{i}$. We need to show that this
relation follows from
$b_{1}b_{2}=b_{3}b_{4}=b_{5}b_{6}$. Cancelling in
$F$, if needed, we may assume that each $b_{i}$
appears at most on one side of the relation. We
also may  assume that not both sides are
divisible in $\langle b_{1},\ldots, b_{6}\rangle$
by one of the equal words $b_{1}b_{2}$,
$b_{3}b_{4}$ and $b_{5}b_{6}$. Further, on both
sides of $v=w$ we need some $b_{i}$ with an even
$i$. Indeed, suppose the contrary,  then $x_{8}$
is not involved in $v$ and $w$. Hence also
$b_{5}$ cannot occur because of $x_{7}$. But, as
$b_{1}$ and $b_{3}$  generate a free abelian
monoid of rank $2$, it then follows that $v$ and
$w$ are identical words in the $b_{i}$'s, as
desired. Hence, by symmetry we may assume that
$w$ contains $b_{2}^{i_{2}}$ with $i_{2}>0$ and
$w$ does not contain $b_{4}$ nor $b_{6}$ as a
factor, and $v$ does not contain $b_{2}$ as a
factor and  contains $b_{4}^{i_{4}}b_{6}^{i_{6}}$
for some nonnegative $i_{4},i_{6}$ with
$i_{2}=i_{4}+i_{6}$ (the latter follows by taking
into the account the degree of $x_{8}$ in the
respective words). Looking at $x_{4}$ we then get
that $b_{1}$ appears in $w$. If $i_{6}=0$ then
(looking at $x_{6}$) we get that $b_{3}$ is in
$v$, a contradiction because
$b_{1}b_{2}=b_{3}b_{4}$ divides then both $v$ and
$w$. So $i_{6}>0$. Also $v$ must contain
$x_{1},x_{5}$ and so $b_{3}$ or $b_{5}$ is in
$v$. The latter is not possible because then $v$
and $w$ are divisible by $b_{1}b_{2}=b_{5}b_{6}$,
a contradiction. Thus, $b_{3}$ occurs in $v$.
Then we must have $i_{4}=0$ because otherwise
$b_{1}b_{2}=b_{3}b_{4}$ divides $v$ and $w$. So
$w=b_{1}^{i_{1}}b_{2}^{i_{2}}b_{5}^{i_{5}}=b_{3}^{i_{3}}b_{6}^{i_{6}}=v$
for some $i_{1},i_{2},i_{6}>0$ and
$i_{3},i_{5}\geq 0$. Then the exponents of
$x_{7}$ show that $i_{2}=i_{5}=0$, a
contradiction. The claim follows. Hence we may
indeed identify $a_{i}$ with $b_{i}$ and $B$ with
$\langle b_{1},\ldots, b_{6}\rangle$.

Next, suppose that
$w=a_{2}^{i_{2}}a_{3}^{i_{3}}a_{4}^{i_{4}}a_{5}^{i_{5}}\in
F\cap \gr (B)$. This is equivalent to the
conditions: $i_{3}+i_{5}\geq 0, i_{3}\geq 0,
i_{4}+i_{5}\geq 0,i_{4}\geq
0,i_{2}+i_{3}+i_{5}\geq 0, i_{2}+i_{3}\geq 0,
i_{2}+i_{4}+i_{5}\geq 0, i_{2}+i_{4}\geq 0$. Let
$j=\min \{i_{3},i_{4}\} \geq 0$. Then
$j+i_{2}+i_{5}\geq 0$. We choose $s,t\geq 0$ so
that $i_{2}+s,i_{5}+t \geq 0$ and $s+t=j$. Then
$w=a_{1}^{s}a_{2}^{i_{2}+s}a_{3}^{i_{3}-j}a_{4}^{i_{4}-j}a_{5}^{i_{5}+t}a_{6}^{t}$.
It is thus clear that $w\in B$. So $B=F\cap \gr
(B)$. Since $F$ is a maximal order it then easily
follows that $B$ is a maximal order.
\end{proof}

We conclude with the promised illustration of Theorem~\ref{thm33}.
This example also shows that Theorem~\ref{thm-primes} cannot be
extended to prime ideals of height exceeding $1$.

\begin{example} \label{mainexample}
Let $K$ be any field and let $R=K\langle x_{1},x_{2},x_{3},x_{4} \rangle $
 be the algebra defined by the following relations:
\begin{eqnarray*} x_{1}x_{4}=x_{2}x_{3},
x_{1}x_{3}=x_{2}x_{4}, x_{3}x_{1}=x_{4}x_{2} \\
x_{3}x_{2}=x_{4}x_{1}, x_{1}x_{2}=x_{3}x_{4},
x_{2}x_{1}=x_{4}x_{3} . \end{eqnarray*} Clearly, $R=K[S]$ for the
monoid $S$ defined by the same presentation. Then $S$ is
cancellative (but the group $SS^{-1}$ is not torsion free) and $R$
is a prime Noetherian PI-algebra that is a maximal order.
Furthermore, there exists a prime ideal $P$ of $S$ so that $K[P]$
is not a prime ideal of $K[S]$.
\end{example}
\begin{proof}
Notice that each of the permutations $(12)(34),
(13)(24)$ and $(14)(23)$ determines an
automorphism of $S$. First we list some
equalities in $S$, namely all relations between
the elements of length $3$. For brevity, we use
the index $i$ in place of the generator $x_{i}$.

 \begin{eqnarray*}
 112=134=244, &
 113=124=344, &
 114=123=343=321=411 \\
  221=243=133, &
  224=213=433, &
  223=214=434=412=322  \\
  331=342=122, &
  334=312=422, &
  332=341=121=143=233 \\
  442=431=211, &
  443=421=311, &
  441=432=212=234=144
\end{eqnarray*}
  \begin{eqnarray*}
   131=142=232=241, &
    141=132=242=231 \\
     313=423=414=324, &
      323=314=424=413 . \end{eqnarray*}
It follows easily that $A=\langle
x_{1}^{2},x_{2}^{2},x_{3}^{2},x_{4}^{2}\rangle $
is an abelian submonoid and it is normal, that is
$xA=Ax$ for every $x\in S$. Moreover
$x_{1}^{2}x_{4}^{2}=x_{2}^{2}x_{3}^{2}$ because
$x_{2}x_{2}x_{3}x_{3}=x_{2}x_{3}x_{3}x_{2}=x_{1}x_{4}x_{4}x_{1}=
x_{4}x_{4}x_{1}x_{1}$.

Let $a_{1}=x_{1}x_{4}, a_{2}=x_{4}x_{1}, a_{3}= x_{1}^{2},
a_{4}=x_{4}^{2}, a_{5}=x_{2}^{2},a_{6}= x_{3}^{2}$ and $B=\langle
a_{1}, a_{2}, a_{3},a_{4},a_{5}, a_{6}\rangle $. From the above
equalities it follows that $B$ is abelian and $sB=Bs$ for every
$s\in S$.

Every element of $S$ that is a word of length $3$
is either of the form $xyy$ or of the form
$xx_{1}x_{4}, xx_{4}x_{1}$ for some $x$. It is
also easy to see that every element of $S$ is of
the form $zu_{1}w_{1}\cdots u_{r}w_{r}$ or
$zu_{1}w_{1}\cdots u_{r}w_{r}u_{r+1}$ or $zc$,
where $z\in A$ and $c\in S$ is an element of
length at most $2$ in the $x_{j}$ and either
$u_{i}\in \{ x_{1},x_{2}\}, w_{i} \in \{
x_{3},x_{4}\} $ or $u_{i}\in \{ x_{3},x_{4}\},
w_{i} \in \{ x_{1},x_{2} \} $ (for all $i$). So
we may assume that $c\in
\{x_{1}x_{2},x_{2}x_{1}\}$. Here $r$ is a
non-negative integer. Using the relations listed
above (especially
$x_{1}x_{4}x_{1}=x_{2}x_{4}x_{2},
x_{4}x_{1}x_{4}=x_{3}x_{1}x_{3}$) it may be
checked that the even powers of
$x_{1}x_{4}=x_{2}x_{3}$ and
$x_{2}x_{4}=x_{1}x_{3}$ are equal and the even
powers of $x_{4}x_{1}=x_{3}x_{2}$ and
$x_{3}x_{1}=x_{4}x_{2}$ are equal. Also
$\{x_{4}x_{1},x_{1}x_{4}\}\{x_{1}x_{2},x_{2}x_{1}\}\subseteq
A\{x_{1}x_{3},x_{3}x_{1}\}$. It follows that
possible forms of elements of $S$ are
\begin{eqnarray}  z(x_{1}x_{4})^{i}, \ z(x_{1}x_{4})^{i}x_{1}, \
z(x_{1}x_{4})^{i}x_{2}, \
z(x_{1}x_{4})^{i}x_{1}x_{3}, \ zx_{1}x_{2}
\label{two}\\ z(x_{4}x_{1})^{i}, \
z(x_{4}x_{1})^{i}x_{4}, \ z(x_{4}x_{1})^{i}x_{3},
\ z(x_{4}x_{1})^{i}x_{3}x_{1}, \ zx_{2}x_{1},
\label{three}
\end{eqnarray}
with $z\in A, i\geq 0$. This leads to
\begin{eqnarray}
S&=&B\cup Bx_{1} \cup Bx_{2}\cup Bx_{3}\cup
Bx_{4} \cup Bx_{1}x_{3}\cup Bx_{3}x_{1} \cup
Bx_{1}x_{2} \cup Bx_{2}x_{1} . \label{semgrpel}
\end{eqnarray}
It follows that $K[S]$ is finite module over the
finitely generated commutative algebra $K[B]$.
Hence $K[S]$ is a Noetherian PI-algebra.

Let $C$ be the free abelian group of rank $4$
generated by elements $a,b,c,d$. Let $X$ be the
free monoid on $x_{1},x_{2},x_{3},x_{4}$.
Consider the monoid homomorphism $\phi
:X\longrightarrow M_{4}(K[C])$ defined by
\begin{eqnarray*}
x_{1}\mapsto  \left(
\begin{array}{cccc}
   0  & a & 0 & 0 \\
   1  & 0 & 0 & 0 \\
   0  & 0 & 0 & a^{-1}bc \\
   0  & 0 & 1 & 0
      \end{array} \right) , &&x_{2}\mapsto \left( \begin{array}{cccc}
    0 & 0 & b & 0 \\
    0 & 0 & 0 & a^{-1}bc  \\
    1 & 0 & 0 &  0 \\
    0 & ab^{-1} & 0 & 0
      \end{array} \right) ,\\
      &&\\
x_{3}\mapsto \left( \begin{array}{cccc}
    0 & 0 & bcd^{-1} & 0 \\
    0 & 0 & 0 & a^{-1}bd \\
    b^{-1}d & 0 & 0 & 0 \\
    0 & ad^{-1} & 0 & 0
      \end{array} \right) ,&&  x_{4}\mapsto \left(
\begin{array}{cccc}
    0 & bcd^{-1} & 0 & 0 \\
    a^{-1}d & 0 & 0 &  0 \\
    0 &  0 & 0 &  d \\
    0 & 0 & ad^{-1}  & 0
      \end{array} \right) .
\end{eqnarray*}
It is easy to check that these matrices satisfy
the defining relations of $S$, so $\phi $ can be
viewed as a homomorphism from $S$ to the group of
monomial matrices over $C$. Moreover
\begin{eqnarray*}a_{3}\mapsto \left(
\begin{array}{cccc}
   a  & 0 & 0 & 0 \\
   0  & a & 0 & 0 \\
   0  & 0 & a^{-1}bc & 0 \\
   0  & 0 & 0 & a^{-1}bc
      \end{array} \right) ,&&
a_{5}\mapsto \left( \begin{array}{cccc}
   b  & 0 & 0 & 0 \\
   0  & c & 0 & 0 \\
   0  & 0 & b & 0 \\
   0  & 0 & 0 & c
      \end{array} \right) \\
      &&\\
a_{4}\mapsto \left( \begin{array}{cccc}
   a^{-1}bc  & 0 & 0 & 0 \\
   0  & a^{-1}bc & 0 & 0 \\
   0  & 0 & a & 0 \\
   0  & 0 & 0 & a
      \end{array} \right) ,&&
a_{6}\mapsto \left( \begin{array}{cccc}
   c  & 0 & 0 & 0 \\
   0  & b & 0 & 0 \\
   0  & 0 & c & 0 \\
   0  & 0 & 0 & b
      \end{array} \right) ,\\
      &&\\
a_{1}\mapsto \left( \begin{array}{cccc}
   d  & 0 & 0 & 0 \\
   0  & bcd^{-1} & 0 & 0 \\
   0  & 0 & bcd^{-1} & 0 \\
   0  & 0 & 0 & d
      \end{array} \right) ,&&
      a_{2}\mapsto \left( \begin{array}{cccc}
   bcd^{-1}  & 0 & 0 & 0 \\
   0  & d & 0 & 0 \\
   0  & 0 & d & 0 \\
   0  & 0 & 0 & bcd^{-1}
      \end{array} \right) .
\end{eqnarray*}
The projection of the group $C'$ generated by $\phi (a_{1}), \phi
(a_{3}), \phi (a_{5}), \phi (a_{6})$ onto the $(1,1)$-entry
contains $a,b,c,d$, whence it is free abelian of rank $4$. So
$C'\cong C$. In particular, $\phi$ is injective on $B_{0}=\langle
a_{1},a_{3},a_{5},a_{6}\rangle$, and thus $B_{0}$ is a free
abelian monoid of rank $4$. Let $B'$ be the abelian monoid with
presentation $\langle y_{1},y_{2},y_{3},y_{4},y_{5},y_{6}\mid
y_{1}y_{2}=y_{3}y_{4}=y_{5}y_{6}\rangle$. Clearly, we have natural
homomorphisms $B'\longrightarrow B\longrightarrow \phi (B)$. We
also know that $\phi (B)$ generates a free abelian group of rank
$\geq 4$ and $B'$ has a group of quotients that is free abelian of
rank $4$. So $B'$ and $\phi(B)$ must be isomorphic, since
otherwise under the map $B'\longrightarrow \phi(B)$ we have to
factor out an additional relation and the rank would decrease. It
follows that $\phi $ is injective on $B$ and thus, because of
Example~\ref{abel-max}, $B$ is a cancellative maximal order with
group of quotients $N= BB^{-1}\cong {\mathbb Z}^{4}$. Note that
$AA^{-1}=\gr (a_{3},a_{5},a_{6})$. Using the defining relations
$a_{1}a_{2}=a_{3}a_{4}=a_{5}a_{6}$ for $B$, it is readily verified
that if $i,j,k\in {\mathbb Z}$ then
$a_{3}^{i}a_{5}^{j}a_{6}^{k}\in A=\langle
a_{3},a_{4},a_{5},a_{6}\rangle$ if and only if $j,k\geq 0$ and
$\min (j,k)+i\geq 0$.

The images under $\phi $ of the first $4$ types
listed in (\ref{two}) have different patterns of
nonzero entries in $M_{4}(K[C])$. The same
applies to the first $4$ types in (\ref{three}).
Notice that $\phi $ is injective on each of the
$10$ types. Suppose that
$s=z(x_{1}x_{4})^{i}x_{1},
t=z'(x_{4}x_{1})^{j}x_{4}$ have the same image.
Then $sx_{4}=z(x_{1}x_{4})^{i+1}$ and $tx_{4}=z
'(x_{4}x_{1})^{j}x_{4}^{2}$ have equal images.
This is not possible because
$(x_{1}x_{4})^{i+1}A\cap
(x_{4}x_{1})^{j}A=\emptyset$ by the above
description of the group $BB^{-1}\cong
\phi(B)\phi(B)^{-1}$. Similarly one deals with
elements of any two different types listed in
(\ref{two}) and (\ref{three}), showing that only
elements the form $s=zx_{1}x_{2},
t=z'x_{2}x_{1}$, where $z,z'\in A$, can have
equal images. Then
$\phi(zx_{1}x_{1}x_{4}x_{1})=\phi
(tx_{4}x_{2})=\phi(sx_{4}x_{2})=\phi(z'x_{4}x_{1}x_{2}x_{2})$
and cancellativity of $\phi (B)\cong B$ yields
$zx_{1}x_{1}=z'x_{2}x_{2}$. Write
$z=a_{3}^{i}a_{5}^{j}a_{6}^{k}$, with $i,j,k\in
{\mathbb Z}$. Since $z,z'\in A$, the above yields
that $z'=a_{3}^{i+1}a_{5}^{j-1}a_{6}^{k}$ and
$j,k\geq 0,-i$ and $j-1,k\geq 0,-(i+1)$. Notice
that
$a_{3}x_{2}x_{1}=x_{1}x_{1}x_{2}x_{1}=x_{2}x_{2}x_{1}x_{2}=a_{5}x_{1}x_{2}$.
Therefore, if $j-1\geq -i$ then
$a_{3}^{i}a_{5}^{j-1}a_{6}^{k}\in A$ and hence
$zx_{1}x_{2}=a_{3}^{i}a_{5}^{j-1}a_{6}^{k}a_{5}x_{1}x_{2}=
a_{3}^{i}a_{5}^{j-1}a_{6}^{k}a_{3}x_{2}x_{1}=z'x_{2}x_{1}$.
On the other hand, if $j=-i$ then
$a_{4}^{j-1}a_{6}^{k-j}\in A$ and
$a_{4}x_{1}x_{2}=a_{6}x_{2}x_{1}$. Hence we also
get
$zx_{1}x_{2}=a_{4}^{j}a_{6}^{k-j}x_{1}x_{2}=a_{4}^{j-1}a_{6}^{k-j}a_{4}x_{1}x_{2}=
a_{4}^{j-1}a_{6}^{k-j}a_{6}x_{2}x_{1}=z'x_{2}x_{1}$.

It follows that $\phi $ is injective on all
elements of types (\ref{two}),(\ref{three}).
Therefore $\phi $ is an embedding and thus $S$ is
cancellative.

We identify $S$ with $\phi (S)$. Put $G=SS^{-1}$. Then $G= N\cup
x_{1} N\cup x_{2}N\cup x_{1}x_{2}N$. Moreover $N\cong {\mathbb
Z}^{4}$ and $N$ is a normal subgroup with $G/N$ the four group.
From (\ref{semgrpel}) it follows that $S\cap N=B$.

We now show that $G$ is dihedral free. To prove
this, notice that $S$ acts by conjugation on $B$
and the generators $x_{i}$ of $S$ correspond to
the following permutations $\sigma _{i}$ of the
generating set of $B$ (the numbers $1,2,3,4,5,6$
correspond to the generators
$a_{1},a_{2},a_{3},a_{4},a_{5},a_{6}$).
\begin{eqnarray} \label{permut}
\sigma_{1} =(12)(56), \, \sigma_{2}= (12)(34),\,
\sigma_{3}= (12)(34),\, \sigma_{4} =(12)(56) .
\end{eqnarray}
Suppose $D\subseteq G$ is an infinite dihedral group such that the normalizer $N_{G}(D)$ of $D$ in
$G$ is of finite index. Let $t\in D$ be an element of order $2$. Then there exists $k\geq 1$ such
that $a_{i}^{k}ta_{i}^{-k}t\in D$ for $i=1,2,\ldots ,6$. Clearly
$$a_{i}^{k}ta_{i}^{-k}t=a_{i}^{k}a_{\sigma(i)}^{k}t^{2}=a_{i}^{k}a_{\sigma(i)}^{k},$$
where $\sigma$ is the automorphism of $N$
determined by $t$. Since $t\in Nx_{1}\cup
Nx_{2}\cup Nx_{1}x_{2}$, it follows that $\sigma
$ is determined by the conjugation by
$x_{1},x_{2}$ or $x_{1}x_{2}$. So $\sigma $
permutes exactly two of the pairs $a_{1},a_{2};
a_{3},a_{4}$ and $a_{5},a_{6}$. It follows that
$a_{i}^{k}a_{j}^{k}, a_{p}^{k}a_{q}^{k}\in D$ for
two different pairs $i,j$ and $p,q$. Hence $\rk
(N\cap D)\geq 2$, a contradiction. Therefore $G$
indeed is a dihedral free group.

Next we show that $\Delta^{+}(G)$ is trivial. For this, let $F$ be
a finite normal subgroup of $G$. Since $N$ is torsion free, it is
clear that $F$ is isomorphic with a subgroup of ${\mathbb Z}
_{2}\times {\mathbb Z} _{2}$. A nontrivial element $t\in F$ must
be of order $2$, whence as above we get that
$a_{i}^{n}a_{j}^{n}t\in F$ for some $i\neq j$ and every $n\geq 1$.
Therefore $F$ is infinite, a contradiction. It follows that
$\Delta ^{+}(G)$ is trivial. Therefore $K[G]$ is prime and hence
$K[S]$ also is prime.

Note that in $G$ we have
$x_{2}^{-1}x_{1}=x_{3}x_{4}^{-1}=x_{4}x_{3}^{-1}$
and
$x_{1}x_{2}^{-1}=x_{3}^{-1}x_{4}=x_{4}^{-1}x_{3}$.
So $x_{2}^{-1}x_{1}$ is an element of order $2$.

We describe the minimal prime ideals of $S$. For
this we first notice that it is easy to see that
the minimal prime ideals of $B$ are:
\begin{eqnarray*}
Q=Q_{1}=(a_{1},a_{3},a_{5}), &&
Q_{2}=Q^{x_{1}}=(a_{2},a_{3},a_{6}),\\
Q_{3}=Q^{x_{2}}=(a_{2},a_{4},a_{5}), &&
Q_{4}=Q^{x_{1}x_{2}}=(a_{1},a_{4},a_{6}),\\
Q'=Q_{5}=(a_{2},a_{3},a_{5}), &&
Q_{6}=(Q')^{x_{1}}=(a_{1},a_{3},a_{6}),\\
Q_{7}=(Q')^{x_{2}}=(a_{1},a_{4},a_{5}),&&
Q_{8}=(Q')^{x_{1}x_{2}}=(a_{2},a_{4},a_{6}).
\end{eqnarray*}
Because of (\ref{permut}), it is easily verified
that $a_{1}a_{2}$ is a central element of $S$ and
that every ideal of $S$ contains a positive power
of $a_{1}a_{2}$. In particular, $a_{1}a_{2}$
belongs to every prime ideal of $S$. Consider in
$B$ the following $G$-invariant ideals:
$$M=(a_{1}a_{3}a_{5},
a_{2}a_{3}a_{6},a_{2}a_{4}a_{5},a_{1}a_{4}a_{6},a_{1}a_{2})$$
and $$M'=(a_{2}a_{3}a_{5},
a_{1}a_{4}a_{5},a_{1}a_{3}a_{6},a_{2}a_{4}a_{6},a_{1}a_{2}).$$
Again because of (\ref{permut}), it is easy to
see that $bSb' \subseteq a_{1}a_{2}S$ for every
defining generator $b$ of $M$ and $b'$ of $M'$.
It follows that a prime ideal of $S$ contains $M$
or $M'$.

Notice that $xa_{1}a_{3}a_{5}y\in xy\{
a_{1}a_{3}a_{5},a_{2}a_{4}a_{5},a_{2}a_{3}a_{6},a_{1}a_{4}a_{6}\}$ for every $x,y\in S$. Therefore
$Sa_{1}a_{3}a_{5}S\cap \langle a_{2},a_{3},a_{5}\rangle=\emptyset$ (for example,
$xya_{1}a_{3}a_{5}\not \in \langle a_{2},a_{3},a_{5}\rangle$ because otherwise
$xy=a_{1}^{-1}a_{3}^{-1}a_{5}^{-1}\langle a_{2},a_{3},a_{5}\rangle \cap (S\cap N)$, which is not
possible because $S\cap N=B$ and $N=\gr(a_{1},a_{2},a_{3},a_{5})$ is free abelian of rank $4$). So
there exists a (unique) ideal $P$ of $S$ that is maximal with respect to the property $P\cap
\langle a_{2},a_{3},a_{5}\rangle =\emptyset$. It is easy to see that $P$ is a prime ideal of $S$.
Since $a_{2}a_{3}a_{5}\not\in P$, we get that $M'\not\subseteq P$, whence $M\subseteq P$.

Let $E=a_{2}a_{3}a_{5}\langle a_{2},a_{3},a_{5}\rangle $. Then
\begin{eqnarray} x_{4}\langle
a_{1},a_{3},a_{6}\rangle x_{1}\subseteq E, \,
x_{2}\langle a_{1},a_{4},a_{5}\rangle x_{2}
\subseteq E,  \, x_{1}x_{3} \langle
a_{2},a_{4},a_{6} \rangle x_{3}x_{1}\subseteq E .
\end{eqnarray}
Therefore $P\cap \langle a_{1},a_{3},a_{6}\rangle=\emptyset, P\cap \langle
a_{1},a_{4},a_{5}\rangle=\emptyset, P\cap \langle a_{2},a_{4},a_{6}\rangle=\emptyset$. For every
element $b\in B\setminus M$ the ideal $bB$ intersects one of the sets
$$a_{2}a_{3}a_{5}\langle a_{2},a_{3},a_{5}\rangle ,a_{1}a_{3}a_{6}\langle
a_{1},a_{3},a_{6}\rangle, a_{1}a_{4}a_{5}\langle a_{1},a_{4},a_{5}\rangle, a_{2}a_{4}a_{6} \langle
a_{2},a_{4},a_{6} \rangle.$$ Since these sets do not intersect $P$, we get that $P\cap B=M$. Hence
it is $G$-invariant and, because $M\subseteq Q_{1}\cap Q_{2}\cap Q_{3}\cap Q_{4}$ and all $Q_{i}$
are minimal primes of $B$, Corollary~\ref{primes-prop} yields that $P\cap B=M=Q_{1}\cap Q_{2}\cap
Q_{3}\cap Q_{4}={\mathcal B}(SMS)$ and $P$ is a minimal prime ideal of $S$. A similar argument
shows that there exists an ideal $P'$ of S that is maximal with respect to the property $P' \cap
\langle a_{1},a_{3},a_{5}\rangle = \emptyset$, and it follows that $P'\cap B=M'=Q_{5}\cap
Q_{6}\cap Q_{7}\cap Q_{8}={\mathcal B}(SM'S)$  is the only other minimal prime of $S$.

In order to continue the proof we first show the following claim
on the representation of elements of $B$.

{\it Claim: Presentation}\\
 An element $b
=a_{1}^{\alpha_{1}}a_{3}^{\alpha_{3}}a_{4}^{\alpha_{4}}a_{6}^{\alpha_{6}}$
of $N=BB^{-1}$ is in $B$ (with each $\alpha_{i} \in {\mathbb Z}$)
if and only if $\alpha_{3},\alpha_{4}\geq 0$ and either (i)
$\alpha_{1}\geq 0$ and $\min (\alpha_{3} ,\alpha_{4})+\alpha_{6}
\geq 0$ or (ii) $\alpha_{1}<0$, $\alpha_{3}+\alpha_{1},
\alpha_{4}+\alpha_{1}\geq 0$ and $\min (\alpha_{3} ,\alpha_{4})
+\alpha_{1} +\alpha_{6}\geq 0$, or equivalently, $\min (\alpha_{3}
+\alpha_{1},\alpha_{4} +\alpha_{1})\geq \max(0, - \alpha_{6})$.

That the first condition is sufficient is easily verified. For the
second condition, we rewrite $a_{1}^{\alpha_{1}}$ as
$a_{3}^{\alpha_{1}}a_{4}^{\alpha_{1}}a_{2}^{-\alpha_{1}}$ and the
result follows because of the first part (by interchanging $a_{1}$
with $a_{2}$).

To prove that they are necessary, suppose
$b=a_{1}^{\alpha_{1}}a_{3}^{\alpha_{3}}
a_{4}^{\alpha_{4}}a_{6}^{\alpha_{6}}\in B$, with each $\alpha_{i}
\in {\mathbb Z}$. Because $Ba_{4}\cap \langle
a_{1},a_{3},a_{6}\rangle =\emptyset$ and $Ba_{3}\cap \langle
a_{1},a_{4},a_{6}\rangle =\emptyset$, it follows that
$\alpha_{3},\alpha_{4}\geq 0$. Suppose now that $\alpha_{1}\geq
0$. Then
$$b=a_{1}^{\alpha_{1}} a_{3}^{\alpha_{3}-\min
(\alpha_{3},\alpha_{4}) } a_{4}^{\alpha_{4}-\min
(\alpha_{3},\alpha_{4}) } a_{5}^{\min (\alpha_{3},\alpha_{4})}
a_{6}^{\alpha_{6}+\min (\alpha_{3},\alpha_{4})}$$ and thus
$a_{1}^{\alpha_{1}} a_{3}^{\alpha_{3}-\min
(\alpha_{3},\alpha_{4})} a_{4}^{\alpha_{4}-\min
(\alpha_{3},\alpha_{4})}a_{5}^{\min (\alpha_{3},\alpha_{4})} \in
Ba_{6}^{-(\alpha_{6}+\min (\alpha_{3},\alpha_{4}))}$. Since the
exponent of $a_{3}$ or $a_{4}$ is $0$, this implies that
$\alpha_{6}+\min (\alpha_{3},\alpha_{4}) \geq 0$, as desired.

On the other hand, suppose that $\alpha_{1}<0$. Then $b=a_{2}^{-\alpha_{1}}
a_{3}^{\alpha_{3}+\alpha_{1}} a_{4}^{\alpha_{4}+\alpha_{1}}
a_{6}^{\alpha_{6}}$. The previous case yields that
$\alpha_{3}+\alpha_{1}\geq 0$, $\alpha_{4}+\alpha_{1} \geq 0$ and
$\min (\alpha_{3}+\alpha_{1},\alpha_{4}+\alpha_{1})
+\alpha_{6}\geq 0$, again as desired. This proves the claim {\it
Presentation}.

Next we show that $S$ is a maximal order. First note that $B$ is a
maximal order by Example~\ref{abel-max}. So, because of
Proposition~\ref{notesjan}, it is sufficient to prove that if
$s\in SS^{-1}\setminus S$ then $\langle S,s \rangle \cap N$
strictly contains $B$.

We know that $G=SS^{-1}= N\cup x_{1} N\cup x_{2}N\cup
x_{1}x_{2}N$. If  $s\in N \setminus B$ then clearly  $\langle S,s
\rangle \cap N$ strictly contains $B$. So, there are three cases
to be dealt with: (1) $s = bx_{1}\in G \setminus S$, (2) $s =
bx_{2}\in G\setminus S$, and (3) $s = bx_{1}x_{2}\in G \setminus
S$, where  $b\in N\setminus B$.

Case (1): $s = bx_{1}\in G \setminus S$. Obviously, $bx_{1}x_{4} =
b a_{1}\in \langle S,s \rangle$. If $ba_{1} \in N\setminus B$,
then we are done. So suppose $ba_{1} \in B$ and $b \in N\setminus
B$. We write $b$ in the following form:
$$a_{1}^{\alpha_{1}}a_{3}^{\alpha_{3}}a_{4}^{\alpha_{4}}a_{6}^{\alpha_{6}}.$$
First, suppose that $\alpha_{1} \geq 0$. By the claim {\it Presentation}, from $ba_{1} \in B$ we
get that $\min(\alpha_{3},\alpha_{4})\geq \max(- \alpha_{6},0)$. But, also from $b\in N\setminus
B$ we get that $\min(\alpha_{3},\alpha_{4}) < \max(- \alpha_{6},0)$, a contradiction. So,
$\alpha_{1}$ has to be strictly negative.

Therefore, assume that $\alpha_{1} < 0$. Then $$b
=a_{1}^{\alpha_{1}}a_{3}^{\alpha_{3}}a_{4}^{\alpha_{4}}a_{6}^{\alpha_{6}}
=
a_{2}^{-\alpha_{1}}a_{3}^{\alpha_{3}+\alpha_{1}}a_{4}^{\alpha_{3}+\alpha_{1}}a_{6}^{\alpha_{6}}$$
and $$ba_{1} = a_{2}^{-\alpha_{1} -
1}a_{3}^{\alpha_{3}+\alpha_{1} +
1}a_{4}^{\alpha_{3}+\alpha_{1}+
1}a_{6}^{\alpha_{6}}.$$ By the claim {\it
Presentation} (by interchanging $a_{1}$ and
$a_{2}$), we get that $\min(\alpha_{3} +
\alpha_{1} + 1,\alpha_{4} + \alpha_{1} + 1)\geq
\max(-\alpha_{6},0)$,  but also $\min(\alpha_{3}
+ \alpha_{1},\alpha_{4} + \alpha_{1}) <
\max(-\alpha_{6},0)$. Hence, $\min(\alpha_{3} +
\alpha_{1} + 1,\alpha_{4} + \alpha_{1} + 1) =
\max(-\alpha_{6},0)$. Suppose $\alpha_{6} \geq
0$. Then $\min(\alpha_{3} + \alpha_{1} +
1,\alpha_{4} + \alpha_{1} + 1) = 0$. Therefore,
$b
=a_{1}^{\alpha_{1}}a_{3}^{\alpha_{3}}a_{4}^{\alpha_{4}}a_{6}^{\alpha_{6}}
= a_{1}^{-1}a_{2}^{\alpha_{3}}a_{4}^{\alpha_{4} -
\alpha_{3}}a_{6}^{\alpha_{6}}$ and $\alpha_{4}>
\alpha_{3}$ or $b
=
 a_{1}^{-1}a_{2}^{\alpha_{4}}a_{3}^{\alpha_{3} -
\alpha_{4}}a_{6}^{\alpha_{6}}$ and $\alpha_{3}\geq \alpha_{4}$. In
the first case we get that $bx_{1}\in S$, since
$x_{4}^{-1}x_{1}^{-1}x_{4}x_{4}x_{1} = x_{4}$. So this case gives
a contradiction and is hence impossible. In the second case,
$bx_{1}x_{1} = a_{1}^{-1}a_{2}^{\alpha_{4}}a_{3}^{\alpha_{3} -
\alpha_{4}+1}a_{6}^{\alpha_{6}} \in (\langle S,s \rangle \cap N
)\setminus B$, as desired. If $\alpha_{6} < 0$, $\min(\alpha_{3} +
\alpha_{1} + 1,\alpha_{4} + \alpha_{1} + 1) = -\alpha_{6}$ and
therefore $b =
a_{1}^{\alpha_{1}}a_{3}^{\alpha_{3}}a_{4}^{\alpha_{4}}a_{6}^{\alpha_{6}}
= a_{1} ^{- 1} a_{2}^{-\alpha_{1} - 1}a_{3}^{\alpha_{3} -
\alpha_{4}}a_{5}^{-\alpha_{6}}$ if $\alpha_{4}\leq \alpha_{3}$ or
$b = a_{1} ^{- 1} a_{2}^{-\alpha_{1} - 1}a_{4}^{\alpha_{4} -
\alpha_{3}}a_{5}^{-\alpha_{6}}$ if $\alpha_{3}< \alpha_{4}$. Since
$-\alpha_{1} - 1 \geq 0$ and $-\alpha_{6} > 0$, by interchanging
$a_{5}$ and $a_{6}$, this case is completely similar to the case
where $\alpha_{6}\geq 0$. This finishes the proof of Case (1).

Case (2): $s=bx_{2}\in G \setminus S$. The permutation $\sigma
=(12)(34)$ determines an automorphism  $\sigma$  on $S$ with
$\sigma (B)=B$. Clearly, $\sigma (s) =b'x_{1}$ with $b'=\sigma
(b)\not\in B$. From Case (1) we get that $\langle S,\sigma
(s)\rangle \cap N $ properly contains $B$. Again applying $\sigma$
to the latter we get that $\langle S,s\rangle \cap N$ properly
contains $B$, as required.

 Case (3):  $s = bx_{1}x_{2}\in G\setminus S$.
Clearly, $bx_{1}x_{2}x_{2}x_{4} = b a_{1}a_{6}\in \langle S,s
\rangle$. If $ba_{1}a_{6} \in N\setminus B$, then we are done. So
suppose $ba_{1}a_{6} \in B$ and $b \in N\setminus B$. Write
$b=a_{1}^{\alpha_{1}}a_{3}^{\alpha_{3}}a_{4}^{\alpha_{4}}a_{6}^{\alpha_{6}}$.
We know that $\alpha_{3},\alpha_{4}\geq0$ and we consider again the two cases:
$\alpha_{1}\geq 0$ and $\alpha_{1} < 0$, separately.

First assume that $\alpha_{1}\geq 0$. Since $ba_{1}a_{6}\in B$ and $b\not \in B$, we get that
$\min(\alpha_{3},\alpha_{4}) = -\alpha_{6} - 1$. Then $b =
a_{1}^{\alpha_{1}}a_{3}^{\alpha_{3}}a_{4}^{\alpha_{4}}a_{6}^{\alpha_{6}} =
a_{1}^{\alpha_{1}}a_{3}^{\alpha_{3} - \alpha_{4}}a_{5}^{-\alpha_{6} - 1}a_{6}^{-1}$, if
$\alpha_{4}\leq \alpha_{3}$, or $b = a_{1}^{\alpha_{1}}a_{4}^{\alpha_{4} -
\alpha_{3}}a_{5}^{-\alpha_{6} - 1}a_{6}^{-1}$, if $\alpha_{3}< \alpha_{4}$. Now, let $\alpha_{1} >
0$. Then $bx_{1}x_{2}\in S$, since $x_{1}x_{4}x_{3}^{-1}x_{3}^{-1}x_{1}x_{2} = x_{2}x_{4}$, a
contradiction. So this case is impossible. Therefore, $\alpha_{1} = 0$. If $\alpha_{3} <
\alpha_{4}$, then $bx_{1}x_{2} \in S$, since $x_{4}x_{4}x_{3}^{-1}x_{3}^{-1}x_{3}x_{4} =
x_{2}x_{1}$, a contradiction, so this case is again impossible. If $\alpha_{4} \leq \alpha_{3}$,
then $bx_{1}x_{2}x_{3}x_{1} = ba_{2}a_{3} \in (\langle S,s \rangle \cap N)\setminus B$, as
desired.

Finally, suppose $\alpha_{1}  < 0$. Then $b
=a_{1}^{\alpha_{1}}a_{3}^{\alpha_{3}}a_{4}^{\alpha_{4}}a_{6}^{\alpha_{6}}
=
a_{2}^{-\alpha_{1}}a_{3}^{\alpha_{3}+\alpha_{1}}a_{4}^{\alpha_{4}+\alpha_{1}}a_{6}^{\alpha_{6}}$
and $ba_{1}a_{6} = a_{2}^{-\alpha_{1} -
1}a_{3}^{\alpha_{3}+\alpha_{1} + 1}a_{4}^{\alpha_{4}+\alpha_{1}+
1}a_{6}^{\alpha_{6}+1}$. By the claim {\it Presentation} (by
interchanging $a_{1}$ and $a_{2}$), we get that $\min(\alpha_{3} +
\alpha_{1} + 1,\alpha_{4} + \alpha_{1} + 1)\geq \max(-\alpha_{6} -
1,0)$,  but also $\min(\alpha_{3} + \alpha_{1},\alpha_{4} +
\alpha_{1}) < \max(-\alpha_{6},0)$. If $\alpha_{6} \geq 0$,
$\min(\alpha_{3},\alpha_{4}) + \alpha_{1} = -1$ and $b =
a_{1}^{-1}a_{2}^{\alpha_{4}}a_{3}^{\alpha_{3} -
\alpha_{4}}a_{6}^{\alpha_{6}}$, if $\alpha_{4}\leq \alpha_{3}$, or
$b = a_{1}^{-1}a_{2}^{\alpha_{3}}a_{4}^{\alpha_{4} -
\alpha_{3}}a_{6}^{\alpha_{6}}$, if $\alpha_{3}<\alpha_{4}$. So, we
have the same conditions here as in Case 1 and the result follows.
If $\alpha_{6} < 0$, there are two possibilities:
$\min(\alpha_{3},\alpha_{4}) + \alpha_{1} = -\alpha_{6} - 1$ or
$\min(\alpha_{3},\alpha_{4}) + \alpha_{1}= -\alpha_{6} - 2$. In
the first case, we get that $b = a_{1}^{-1}a_{2}^{-\alpha_{1} -
1}a_{3}^{\alpha_{3} - \alpha_{4}}a_{5}^{-\alpha_{6}}$ if
$\alpha_{4}\leq \alpha_{3}$, or $b = a_{1}^{-1}a_{2}^{-\alpha_{1}
- 1}a_{4}^{\alpha_{4} - \alpha_{3}}a_{5}^{-\alpha_{6}}$ if
$\alpha_{3}<\alpha_{4}$. Since $-\alpha_{1} -1 \geq 0$ and $-
\alpha_{6}> 0$, by interchanging $a_{5}$ and $a_{6}$, this case is
completely similar to the case where $\alpha_{6}\geq 0$ and hence
can also be treated as in Case 1. Finally, if the second
possibility holds, that is $\min(\alpha_{3},\alpha_{4}) +
\alpha_{1}= -\alpha_{6} - 2$, then $b = a_{1}^{- 2}a_{2}^{- 2 -
\alpha_{1}}a_{3}^{\alpha_{3} -\alpha_{4}}a_{5}^{-\alpha_{6}}$, if
$\alpha_{4}\leq \alpha_{3}$, or $b = a_{1}^{- 2}a_{2}^{- 2 -
\alpha_{1}}a_{4}^{\alpha_{4} -\alpha_{3}}a_{5}^{-\alpha_{6}}$, if
$\alpha_{4} > \alpha_{3}$, so we always get that $ba_{2}a_{3} \in
(\langle S,s \rangle \cap N)\setminus B$, as desired.

This finishes the proof of the fact that $S$ is a maximal order.
Since $P\cap B$ is invariant for every minimal prime $P$ of $S$,
it then follows from Theorem~\ref{thm33} that $K[S]$ is a maximal
order.

Let $V$ be the ideal of $S$ generated by the elements $a_{3},a_{4},a_{5},a_{6}, x_{1}x_{2},
x_{2}x_{1}$. It easily follows that the elements of $S\setminus V$ are of the form
\begin{eqnarray*}  (x_{1}x_{4})^{i}, \ (x_{1}x_{4})^{i}x_{1}, \
(x_{1}x_{4})^{i}x_{2}, \
(x_{1}x_{4})^{i}x_{1}x_{3}, \\
(x_{4}x_{1})^{i}, \ (x_{4}x_{1})^{i}x_{4}, \ (x_{4}x_{1})^{i}x_{3}, \ (x_{4}x_{1})^{i}x_{3}x_{1},
\end{eqnarray*}
where $i$ is a non-negative integer. Then $S\setminus V= \{ 1\}\cup I$ where the set $I$ can be
written in matrix format as a union of disjoint sets:
 $$\left( \begin{array}{cc}
    I_{11} & I_{12}\\
    I_{21} & I_{22}
 \end{array} \right) ,$$
with $I_{11} =\langle a_{2}\rangle a_{2}\cup \langle a_{2}\rangle
x_{3}x_{1}$, $I_{12}=\langle a_{2}\rangle x_{3} \cup \langle
a_{2}\rangle x_{4}$, $I_{21} = x_{1} \langle a_{2}\rangle \cup
x_{2}\langle a_{2}\rangle$ and $I_{22} =\langle a_{1}\rangle
a_{1}\cup \langle a_{1}\rangle x_{1}x_{3}$. Moreover,
$I_{ij}I_{kl}\subseteq I_{il}$ if $j=k$, and is contained in $V$
otherwise. In $S/V$ the set $I\cup \{ 0\}$ is an ideal and the
semigroup $I_{11}$ (treated as a subsemigroup of $S$) has a group
of quotients $H=\gr (a_{2},x_{3}x_{1})$. Since
$a_{2}^{2}=(x_{3}x_{1})^{2}$ and $a_{2} (x_{3}x_{1})
=(x_{3}x_{1})a_{2}$, we get that $H$ is isomorphic with ${\mathbb
Z} \times {\mathbb Z}_{2}$. From the matrix pattern of $I$ it
follows that $S/V$ is a prime semigroup. So $V$ is a prime ideal
of $S$. However, because $K[H]$ and thus $K[I_{11}]$ is not prime,
standard generalized matrix ring arguments yield that $K[S]/K[V]$
is not prime.
\end{proof}

\vspace{20pt}

 \noindent
 \begin{tabular}{ll}
 I. Goffa and E. Jespers & J. Okni\'{n}ski\\
 Department of Mathematics& Institute of Mathematics\\
 Vrije Universiteit Brussel & Warsaw University\\
 Pleinlaan 2& Banacha 2\\
 1050 Brussel, Belgium& 02-097 Warsaw, Poland\\
 efjesper@vub.ac.be and igoffa@vub.ac.be & okninski@mimuw.edu.pl
 \end{tabular}

\end{document}